\documentclass[11pt]{article}

\usepackage[OT1]{fontenc}
\usepackage{amsthm,amsmath,amssymb}
\usepackage{multirow}
\usepackage{fullpage}
\usepackage{url}
\usepackage{smile}
\usepackage{pifont}
\usepackage{pdflscape}
\setcitestyle{round,semicolon}
\usepackage[colorlinks,
            linkcolor=red,
            anchorcolor=blue,
            citecolor=blue
            ]{hyperref}

\newcommand{\tickYes}{\checkmark}
\newcommand{\tickNo}{\hspace{1pt}\ding{55}}

\DeclareMathOperator*{\argsup}{argsup}
\DeclareMathOperator*{\arginf}{arginf}

\numberwithin{equation}{section}
\numberwithin{theorem}{section}

\begin{document}

\title{\huge A topologically valid definition of depth for functional data}

\author{Alicia Nieto-Reyes\thanks{Departamento de Matem\'aticas, Estad\'istica y Computaci\'on, Universidad de Cantabria; \texttt{alicia.nieto@unican.es}}
%\\ \small{\emph{Universidad de Cantabria}}  
\and Heather Battey\thanks{Current addresses: (1) ORFE, Princeton University; \texttt{hbattey@princeton.edu}; \newline 
\textcolor{white}{white space ....... space} (2) Department of Mathematics, Imperial College London; \texttt{h.battey@imperial.ac.uk}}
%\\ \small{\emph{University of Bristol}}
}
%\date{}

\normalsize{ }

\maketitle
%\doublespacing
\begin{abstract}
The main focus of this work is on providing a formal definition of statistical depth for functional data on the basis of six properties, recognising topological features such as continuity, smoothness and contiguity. Amongst our depth defining properties is one that addresses the delicate challenge of inherent partial observability of functional data, with fulfilment giving rise to a minimal guarantee on the performance of the empirical depth beyond the idealised and practically infeasible case of full observability. As an incidental product, functional depths satisfying our definition achieve a robustness that is commonly ascribed to depth, despite the absence of a formal guarantee in the multivariate definition of depth. We demonstrate the fulfilment or otherwise of our properties for six widely used functional depth proposals, thereby providing a systematic basis for selection of a depth function. % and elucidating \ldots
\end{abstract}

%\hspace{.3cm}\emph{Keywords:}  functional data, multivariate statistics, partial observability, robustness, statistical depth.

\section{Introduction}\label{In}
This work intersects the areas of functional data analysis (FDA) and statistical depth. FDA provides an alternative way of studying traditional data objects, recognising that it is sometimes more natural and more fruitful to view a collection of measurements as partially observed realisations of random functions. Prototypical examples of functional data objects include growth trajectories, handwriting data and brain imaging data. On the other hand, statistical depth (henceforth referred to as depth) is a powerful data analytic and inferential tool, able to reveal diverse features of the underlying distribution such as spread, shape and symmetry \citep{LiuSingh1999}. The ability of depth to reveal distributional features has been exploited in novel ways to define, inter alia, depth-based classifiers \citep[e.g.][]{LiuCA2012,Paindaveine2014}.

The main focus of this work is on providing a formal definition of depth for functional data, justified on the basis of several properties. The definition fills an important void in the existing literature because na\"ive extensions of multivariate depth constructions, designed to satisfy the properties deemed suitable in multivariate space, neglect the topological features of functional data and often give rise to absurd depth computations \citep{Indio,IndioD}. The need for such a definition was first pointed out in the conference proceedings \cite{Iwfos}, where a crude first attempt to address the problem was made. Undesirable behaviour is also evident for specific constructions of functional depth examples that have been proposed  without suitable reflection on the properties sought \citep{Rebecka}.

The properties that constitute our definition, not only provide a sophisticated extension of those defining the multivariate depth, recognising topological features such as continuity, contiguity and smoothness, but also implicitly address several common or inherent difficulties associated with functional data. Amongst our six depth defining properties is one that tackles the delicate challenge of inherent partial observability of functional data, providing a minimal guarantee on the performance of the empirical depth beyond the idealised and practically infeasible case of full observability. Robustness to the presence of outliers is often cited as one of the defining features of empirical depth \citep{Rebecka}. Our definition of functional depth automatically yields a robust estimator of the population depth in the sense of qualitative robustness \citep{HampelAoS}. As we elucidate in Section \ref{subSecDiscussProperties}, none of the properties constituting the multivariate definition of depth \citep{Zuo00, Liu90} give rise to this property, thus a further contribution of our work is the insight that the existing definition for the multivariate framework is insufficient to guarantee robustness of the multivariate empirical depth. A further challenge, automatically addressed (if present) by our definition, pertains to functional data  exhibiting little variability over a subset of the domain and significantly overlapping one another on this set. Intuitively, functional observations over such a domain ought to play a reduced role in the assignment of depth \citep{HubertC}, especially in light of the partial observability and the convention to pre-process the partial observations.

We demonstrate the fulfilment or otherwise of our depth defining properties for six widely used functional depth functions, from which we conclude that the $\emph{$h$-depth}$ \citep{Cuevas07} is the most well-reasoned in terms of number of properties satisfied.

The remainder of this paper is organized as follows. Section \ref{Sbackg} provides an explanation of the notion of depth at the heuristic level, tracking its chronological development, before providing the formal definition of depth in $\mathbb{R}^{p}$, $p\geq 1$, as set forth in \citet{Zuo00, Liu90}. Section \ref{Sbackg} also formalises the functional data setting and defines the notation used in the paper. A formal definition of depth in function space appears in Section \ref{sectionDefinition}, together with a justification of the properties upon which it is based and a thorough discussion of their implications. Section \ref{SectionCompareDepths} analyses existing constructions of functional depth, establishing the fulfilment or otherwise of each property appearing in the definition of functional depth. All the proofs appear in Section \ref{sectionProofs}.

\section{Background and notation}\label{Sbackg}

\subsection{Historical development and a heuristic explanation of depth}

Unlike the univariate case in which there is no ambiguity in the definition of order, when data provide coordinates in a higher dimensional space the notion of order is ill-defined; for instance, in $\mathbb{R}^{2}$ it is not clear whether (3,6) is larger or smaller than (5,4). This fact led to a body of work in the 1970s, proposing new exploratory data analysis tools for assigning ranks to points in a data set. The  method of \emph{convex hull peeling}, credited to J.~W.~Tukey \citep{Huber1972,Barnett1976} is a particularly intuitive example. A pedagogical description of the procedure for the bivariate case is provided in \citet{PeterBook}, where readers are encouraged to envisage the data points as pins on a board. A large elastic band is looped around the pins forming the convex hull of the data points. The data points touching the elastic band are the extremes of the empirical distribution and are assigned rank one and discarded. The procedure is repeated to identify the next most extreme points, which are assigned rank two, and so on. Clearly in this example, the empirical distribution plays an important role in the assignment of rank, where roughly speaking, data points closer to the centre of the empirical distribution receive higher rank(s), giving rise to a centre-outward ordering.

J.~W.~Tukey coined the term \emph{depth} in \cite{Tukey} as the collection of exploratory procedures for assigning ranks to points in a data set. There, he proposed the celebrated \emph{halfspace depth}, or \emph{Tukey depth}, of a data point in $\mathbb{R}^p$ with respect to (henceforth w.r.t.) a multidimensional sample. \cite{Rousseeuw} later defined the halfspace depth w.r.t.~a generic measure as opposed to the empirical measure, broadening the purely data analytic perspective. Thus, modern usage of the term depth refers to a much more general class of objects. The underlying mathematical idea behind these depth constructions and others, is that a probability measure maps events in the Borel $\sigma$-algebra to [0,1], a space on which the assignment of order poses no concern.

Since Tukey's seminal work, many alternative examples of depth have been proposed. It was, however, the \emph{simplicial depth} \citep{Liu90} that sparked a resurgence of research on the topic throughout the 90s and 00s. Simplicial depth is shown in \citet{Liu90} to possess several desirable properties, on the basis of which the definition of depth is formalised in \citet{Zuo00}, reproduced in Definition \ref{zuoDefinition} for ease of reference. In Definition \ref{zuoDefinition}, $\mathcal{P}$ denotes the class of distributions on the Borel sets of $\mathbb{R}^p$, and $P=P_X$ denotes the distribution of a general random vector $X$; the subscript $X$ is suppressed when there is no need to be explicit.

\begin{definition}[\cite{Zuo00, Liu90}] \label{zuoDefinition}
The bounded and non-negative mapping
$D(\cdot,\cdot):\mathbb{R}^p\times\mathcal{P}\longrightarrow \mathbb{R}$ is
called a \emph{statistical depth function} if it satisfies the
following properties:

\begin{enumerate}
\item \emph{Affine invariance.}
$D(Ax+b,P_{AX+b})=D(x,P_X)$ holds for any $\mathbb{R}^p$-valued random vector $X$, any $p\times p$ non-singular matrix $A$ and any $b \in \mathbb{R}^p.$
\item \emph{Maximality at centre.}
$D(\theta,P)=\sup_{x\in\mathbb{R}^p}D(x,P)$ holds for any $P\in\mathcal{P}$
having a unique centre of symmetry $\theta$ w.r.t.~some notion of symmetry.
\item \emph{Monotonicity relative to the deepest point.}
For any $P\in\mathcal{P}$ having deepest point $\theta,$
$D(x,P)\leq D(\theta+\alpha(x-\theta),P)$ holds for all $\alpha\in
[0,1].$
\item \emph{Vanishing at infinity.}
$D(x,P)\rightarrow 0$ as $\| x\| \rightarrow \infty,$ for each $P\in\mathcal{P}$, where $\|\cdot\|$ is the Euclidean norm.
\end{enumerate}
\end{definition}
For a discussion of centre of symmetry in $\mathbb{R}^{p}$, see \cite{Zuo00}; a more general discussion, applicable to function spaces, is provided in Section \ref{discussP2}. Four further properties purported in \cite{SerflingDIMACS} as desirable but not necessary are reproduced in (i)-(iv) below. 

\begin{itemize}
\item[(i)] \emph{Symmetry.} If $P$ is symmetric about $\theta$ in some sense, then so is $D(x,P).$
\item[(ii)] \emph{Continuity of $D(x,P)$ as a function of $x$.} Or merely \emph{upper semi-continuity.}
\item[(iii)] \emph{Continuity of $D(x,P)$ as a function of $P$.}
\item[(iv)] \emph{Quasi-concavity as a function of $x.$} The set $\{x:D(x,P)\geq c\}$ is convex for each real $c.$
\end{itemize}

Upper semicontinuity is a weaker requirement than continuity. In $\RR^{d}$, it is natural to obligate the depth function to preserve the upper semicontinuity property of the distribution function. This statement has a straightforward extension to function spaces, which is addressed in Section \ref{sectionDefinition}. (iii), although not required to provide an order, is indispensable in view of the fact that statisticians do not have access to the true $P$ but rather an empirical counterpart, which converges almost surely to $P$ uniformly over the Borel sets of the domain on which $P$ is defined. It has yet further important implications, explained in Section \ref{discussP6}. Suitable modifications of Properties 1-4, as well as (ii) and (iii) are thus amongst our defining properties of functional depth and are thoroughly justified and discussed in the functional framework in Section \ref{subSecDiscussProperties}. No further attention is dedicated to properties (i) and (iv), which are not deemed necessary, neither in the multivariate nor functional framework. The designation of $D$ as a bounded and non-negative mapping is also unnecessary to provide an order and thus does not appear in our definition of functional depth in Section \ref{sectionDefinition}.

\subsection{The FDA framework}\label{fdaFramework}

To formalise the  FDA framework, a data point is thought of as a realisation of the random function $\{X(v):v\in\mathcal{V}\}$, where $\mathcal{V}$ is a compact subset of $\mathbb{R}^{d}$ for $d\geq 1$. Letting $\Omega$ denote the underlying sample space, $\{X(v): v\in\mathcal{V}\}:=\{X(\omega,v): \omega\in\Omega, v\in\mathcal{V}\}$ is the map $X: \Omega \rightarrow \mathfrak{F}$, where $\mathfrak{F}$ is a function space, whilst for a fixed $\omega \in \Omega$, $X(\omega,\cdot)$ maps from $\mathcal{V}$ to a vector space $\mathbb{F}$. There is a rich body of work concerning $\mathfrak{F}=\mathbb{L}_{2}(\mathcal{V},\lambda)$, the space of Lebesgue square integrable functions from $\mathcal{V}$ to $\mathbb{F}=\mathbb{R}$ (here and henceforth, $\lambda$ denotes Lebesgue measure on $\mathcal{V}$). Non-standard choices of $\mathfrak{F}$ will undoubtedly become more prevalent in the FDA literature, which currently accommodates functional manifolds \citep{functmanifold}, and multivariate functional spaces $\mathfrak{F}=\otimes_{k=1}^{K}\mathbb{L}_{2}(\mathcal{V}_{k},\lambda)$ \citep{ChiouMuller2014} as well as a variety of smoothness classes embedded in $\mathbb{L}_{2}(\mathcal{V},\lambda)$. In the interest of generality, for the definition of functional depth, we do not restrict $\mathfrak{F}$ beyond the assumption that there exists a metric $d$ on $\mathfrak{F}$ such that $(\mathfrak{F},d)$ is a separable metric space.

A further distinguishing feature of functional data is that they are inherently partially observed. Although theoretically infinite dimensional data objects, due to the limitations of the data collection instruments or the experimental design, each functional data object is only ever recorded at a finite set of discretisation points, which we denote by $\mathcal{V}'\subset \mathcal{V}$.

The following notation is henceforth used throughout. $(\mathfrak{F},d)$ is a separable metric space and $\mathcal{A}$ is the $\sigma$-algebra on $\mathfrak{F}$ generated by the open $d$ metric balls. Separability of $(\mathfrak{F},d)$ guarantees that $\mathcal{A}$ coincides with the Borel $\sigma$-algebra on $\mathfrak{F}$ \citep[see e.g.][Chapter 1.7]{vdVW1996}. $(\mathfrak{F},\mathcal{A},P)$ is a probability space with $P\in\mathcal{P}$, the space of all probability measures on the Borel sets of $\mathcal{A}$. Particular instances of $\mathfrak{F}$ to which reference is made are: $\mathbb{H}(\mathcal{V})$, an infinite dimensional Hilbert space on $\mathcal{V}$; $\mathcal{C}(\mathcal{V})$, the space of continuous functions on $\mathcal{V}$; $\mathbb{L}_{p}(\mathcal{V},\lambda)$, the space of Lebesgue $p$-integrable functions on $\mathcal{V}$, where $1\leq p<\infty$; $\mathbb{L}_{\infty}(\mathcal{V})$, the space of uniformly bounded functions on $\mathcal{V}$; and $\mathfrak{W}^{k,p}(\mathcal{V},\lambda)$, the Sobolev space of Lebesgue $p$-integrable functions on $\mathcal{V}$ whose weak derivatives up to order $k\geq 1$ are Lebesgue $p$-integrable on $\mathcal{V}$, where $1\leq p<\infty$. To avoid excessive notation, unless otherwise stated, the argument(s) $\mathcal{V}$ and $\lambda$ (if applicable) are tacit when we write $\mathcal{C}$, $\mathbb{H}$, $\mathbb{L}_{p}$, $\mathbb{L}_{\infty}$ and $\mathfrak{W}^{k,p}$. Similarly, $\|x\|_{\mathbb{L}_{p}(\mathcal{V},\lambda)}=\bigl(\int_{\mathcal{V}} x(v)^{p}\lambda(dv)\bigr)^{1/p}$ is henceforth referred to in the more compact form $\|x\|_{\mathbb{L}_{p}}$. In normed spaces, the metric $d$ will most naturally be a norm; in this case $d=\|\cdot\|_{\mathbb{L}_{p}}$ is used to mean $d(x,y)=\|x-y\|_{\mathbb{L}_{p}}$. $\mathbb{H}$ is most naturally endowed with its inner product norm $\|x-y\|_{\mathbb{L}_{2}}=\sqrt{\langle x-y,x-y\rangle}$ for $x,y\in\mathbb{H}$, whilst $\mathbb{L}_{\infty}$ is most naturally endowed with the supremum norm $\|x-y\|_{\infty}=\sup_{v\in\mathcal{V}}|x(v)-y(v)|$ for $x,y\in\mathbb{L}_{\infty}$. Recall from the above introduction to the FDA framework that for any $\omega\in\Omega$, $X(\omega,\cdot):\mathcal{V}\rightarrow \mathbb{F}$, where $\mathbb{F}$ is a vector space; unless otherwise stated, $\|\cdot\|$ will be used to denote an arbitrary norm on $\mathbb{F}$. For any $x\in\mathfrak{F}$, $x(H):=\{x(v):v\in H\subseteq \mathcal{V}\}$, whilst $x(\mathcal{V})$ is tacitly implied by $x$. Finally, a sample $X_{1},\ldots,X_{n}$ of random draws from $P$ gives rise to the empirical measure $P_{n}$, a collection of $\frac{1}{n}$-weighted point masses at $X_{1},\ldots,X_{n}$. $\widetilde{P}_{n}$ is used to denote the empirical measure of $\widetilde{X}_{1},\ldots,\widetilde{X}_{n}$, which is a sample of reconstructed functional data objects based on the random sample $\{X_{i}(\mathcal{V}'_{i}): i=1,\ldots,n\}$ of partially observed functional data objects, where $\mathcal{V}'_{i}\subset\mathcal{V}$ is a finite set that may be different for every $i\in\{1,\ldots,n\}$.

\section{Formal definition of functional depth}\label{sectionDefinition}

The definition of functional depth provided in this section refers to the concept of centre of symmetry, which is elucidated in Section \ref{discussP2}, and relies on the following preliminary definition.

\begin{definition}\label{definitionConvHull}
Let $(\mathfrak{F},\mathcal{A},P)$ be a probability space as in Section \ref{fdaFramework}. Define $\mathcal{E}$ to be the smallest set in the $\sigma$-algebra $\mathcal{A}$ such that $P(\mathcal{E})=P(\mathfrak{F})$. Then the convex hull of $\mathfrak{F}$ with respect to $P$ is defined as
\[
\mathfrak{C}(\mathfrak{F},P) :=\left\{x\in\mathfrak{F}: x(v)=\alpha L(v)+(1-\alpha)U(v): v\in\mathcal{V}, \alpha\in[0,1]\right\},
\]
where $\displaystyle{U := \{\sup_{x\in\mathcal{E}} x(v): v\in\mathcal{V}\}}$ and $\displaystyle{L := \{\inf_{x\in\mathcal{E}} x(v): v\in\mathcal{V}\}}$.
\end{definition}

\begin{definition}\label{definitionFDepth}
Let $(\mathfrak{F},\mathcal{A},P)$ be a probability space as in Section \ref{fdaFramework}. Let $\mathcal{P}$ be the space of all probability measures on $\mathfrak{F}$. The mapping
$D(\cdot,\cdot):\mathfrak{F}\times\mathcal{P}\longrightarrow \mathbb{R}$ is a \emph{statistical functional depth} if it satisfies
properties P-1. to P-6, below.
\end{definition}

\begin{enumerate}
\item[P-1.]  \emph{Distance invariance}. $D(f(x),P_{f(X)})=D(x,P_{X})$ for any $x\in\mathfrak{F}$ and $f:\mathfrak{F}\rightarrow\mathfrak{F}$ such that for any $y\in\mathfrak{F}$, $d(f(x),f(y))=a_f \cdot d(x,y)$, with $a_f\in\mathbb{R}\backslash\{0\}$.
\end{enumerate}

\begin{enumerate}
\item[P-2.] \emph{Maximality at centre}.
For any $P\in\mathcal{P}$ possessing a unique centre of symmetry  $\theta\in\mathfrak{F}$ w.r.t.~some notion of functional symmetry, $D(\theta,P)=\sup_{x\in\mathfrak{F}}D(x,P)$.
\end{enumerate}

\begin{enumerate}
\item[P-3.] \emph{Strictly decreasing with respect to the deepest point}. 
For any $P\in\mathcal{P}$ such that $D(z,P)=\max_{x\in\mathfrak{F}}D(x,P)$ exists, $D(x,P)< D(y,P)<D(z,P)$ holds for any $x,y\in\mathfrak{F}$  such that $\min\{d(y,z),d(y,x)\}>0$ and  $\max\{d(y,z),d(y,x)\}<d(x,z).$
\end{enumerate}

\begin{enumerate}
\item[P-4.] \emph{Upper semi-continuity in $x$}. $D(x,P)$ is upper semi-continuous as a function of $x,$ i.e., for all $x\in\mathfrak{F}$ and for all $\epsilon>0$, there exists a $\delta>0$ such that 
\begin{equation}\label{checkUSC}
\sup_{y: d(x,y)<\delta}D(y,P)\leq D(x,P)+\epsilon.
\end{equation}
\end{enumerate}

\begin{enumerate}
\item[P-5.] \emph{Receptivity to convex hull width across the domain.} $D(x,P_X)<D(f(x),P_{f(X)})$ for any $x\in \mathfrak{C}(\mathfrak{F},P)$ with $D(x,P)<\sup_{y\in\mathfrak{F}}D(y,P)$ and $f:\mathfrak{F}\rightarrow\mathfrak{F}$ such that $f(y(v))=\alpha(v) y(v)$ with $\alpha(v)\in (0,1)$ for all $v\in L_\delta$ and $\alpha(v)=1$ for all $v\in L^c_\delta.$ 
\[
L_{\delta}:=\argsup_{H\subseteq\mathcal{V}} \Bigl\{\sup_{x,y\in\mathfrak{C}(\mathfrak{F},P)} d(x(H),y(H))\leq\delta \Bigr\}
\]
for any $\delta\in[\inf_{v\in\mathcal{V}} d(L(v),U(v)),d(L,U))$ such that $\lambda(L_{\delta})>0$ and $\lambda(L^c_{\delta})>0.$
\end{enumerate}

\begin{enumerate}
\item[P-6.] \emph{Continuity in $P$}. For all $x\in\mathfrak{F}$, for all $P\in\mathcal{P}$ and for every $\epsilon>0$, there exists a $\delta(\epsilon)>0$ such that $|D(x,Q)-D(x,P)|<\epsilon$ $P$-almost surely for all $Q\in\mathcal{P}$ with $d_{\mathcal{P}}(Q,P)<\delta$ $P$-almost surely, where $d_{\mathcal{P}}$ metricises the topology of weak convergence.
\end{enumerate}

\subsection{Discussion of the functional depth defining properties}\label{subSecDiscussProperties}

\subsubsection{Discussion of P-1.~Distance invariance.}\label{discussP1} Property P-1.~is the generalisation from $\mathbb{R}^{d}$ to $\mathfrak{F}$ of Property 1 of \citet{Zuo00}, also considered in Theorem 3 of \citet{Liu90}. It states that any mapping from $\mathfrak{F}$ to $\mathfrak{F}$ that preserves, up to a scaling factor, the relative distances between elements in the $d$ metric, also preserves the depth in the transformed space. As an example, consider ($\mathfrak{F},d$)$=$($\mathbb{L}_{2},\|\cdot\|_{\mathbb{L}_{2}}$) and suppose $\mu:=\mathbb{E}(X)=\int x P(dx)$ is known. Then Property P-1.~ensures that the depth is unaffected by recentring around the zero function because $\|x-y\|_{\mathbb{L}_{2}}=\|(x-\mu)-(y-\mu)\|_{\mathbb{L}_{2}}$ for all $x,y \in\mathbb{L}_{2}$.

\subsubsection{Discussion of P-2.~Maximality at centre.}\label{discussP2} P-2.~is the most logically contentious of the properties listed. The reason is that, even for distributions on $\mathbb{R}^{d}$, there is no unique notion of symmetry, a fact that is a fortiori true in function spaces. Indeed, since depth itself was originally conceived as a way to give meaning to the concept of centre of symmetry, the deepest element is no less valid as a centre of symmetry than any other definition, giving rise to the somewhat paradoxical conclusion that P-2.~is always achieved with $\theta$ equal to the deepest point, as long as $\sup_{x\in\mathfrak{F}}D(x,P)=\max_{x\in\mathfrak{F}}D(x,P)$. It is more meaningful to consider the behaviour of $D$ for a particular $P$ for which many notions of centre of symmetry coincide at $\theta$. In $\mathbb{R}$ such a $P$ is the Gaussian distribution, for which the median is equal to the mean and is a centre of symmetry with respect to many notions of symmetry including central symmetry and halfspace symmetry \citep[e.g.][]{SerflingS}. In the setting of $\mathfrak{F}=\mathbb{H}$, the analogue of the Gaussian distribution is the Gaussian process. With this in mind, verification of the following property is insightful.

\begin{enumerate}
\item[P-2G.] \emph{Maximality at Gaussian process mean}.
For $P$ a zero-mean, stationary, almost surely continuous Gaussian process on $\mathcal{V}$, $D(\theta,P)=\sup_{x\in\mathfrak{F}}D(x,P)\neq\inf_{x\in\mathfrak{F}}D(x,P)$, where $\theta$ is the zero mean function.
\end{enumerate}

\begin{remark}
Existence of $\mathbb{E}X$ is guaranteed when $X\sim P$ with $P$ a Gaussian process.
\end{remark}

The zero function of Property P2-G is the centre of symmetry of the mean zero Gaussian process with respect to all notions of functional symmetry that have been tacitly introduced via existing depth constructions, for instance pointwise angular symmetry in \citet{Fraiman} and \citet{Romo09}, and pointwise halfspace symmetry in \citet{HubertC}. If a distribution $P_{X}$ on $\mathfrak{F}$ is pointwise halfspace symmetric about $z$, then for every $v\in\mathcal{V}$, the corresponding distribution of $X(v)$ is halfspace symmetric around $z(v)$. 

Property P-2., in partnership with P-3., leads to the centre-outward ordering for which depth was originally conceived. Outward orderings from local centres of symmetry are also possible \citep[see][]{Brusel}, and are induced by constructions that attach greater importance to probabilities $P(A)$ for Borel sets $A$ to which the evaluation point $x$ have close proximity, where proximity is measured by a suitable metric. The relative weighting depends on the features of $P$ that one would like to detect through the use of the local depth function. As the weighting rule becomes close to uniform, the local features are blurred, resulting in global behaviour of any local depth construction. Local centre-outward orderings are not induced by our definition.

\subsubsection{Discussion of P-3.~Strictly decreasing with respect to the deepest point.}\label{discussP3} For some function spaces $\mathfrak{F}$, there is more than one natural metric $d$. For instance, if $\mathfrak{F}=\mathbb{L}_{\infty}\cap\mathfrak{W}^{k,2}$, $(\mathfrak{F},d)$ is separable with respect to the supremum norm, the standard Sobolev inner product norm \citep{Adams1975}, or its slight generalisation, as employed in \citet{Silverman1996}. With this example in mind, setting $d=\|\cdot\|_{\infty}$ and $\mathcal{V}\subset \mathbb{R}$, Property P-3.~ensures that the depth prescribes successively lower depths to functions that only belong to successively larger envelopes around the deepest point $z$. However, when $d$ is the standard Sobolev inner product norm the depth prescribes successively lower depths to functions which lie in successively larger Sobolev balls around $z$, i.e.~its prescription takes account of the distance of $x$ from $z$ in derivative space as well as in $\mathbb{L}_{2}$ norm, assigning low depth to functions much rougher than $z$.

P-3.~has two further implications. The first is that 
\begin{equation}\label{eqLimInf}
\lim_{x:\;d(x,z)\rightarrow \infty}D(x,P)=\inf_{x\in\mathfrak{F}}D(x,P),
\end{equation}
where $z$ is such that $D(z,P)=\max_{x\in\mathfrak{F}} D(x,P)$ exists and where the convention in current literature is to construct $D(\cdot,P)$ such that $\inf_{x\in\mathfrak{F}}D(x,P)=0$ for any $P\in\mathcal{P}$. Equation \eqref{eqLimInf} itself leads to the conclusion of Lemma \ref{lemmaLpSpaces}. 
\begin{lemma}\label{lemmaLpSpaces}
Let $(\mathfrak{F},d)$ be a functional metric space such that $d=\|\cdot\|_{\mathbb{L}_{p}}$, then for each $P\in\mathcal{P}$, \eqref{eqLimInf} implies that $D(x,P)\rightarrow \inf_{x\in\mathfrak{F}}D(x,P)$ as $\|x(v)\|
\rightarrow \infty$ for Lebesgue almost every  $v\in \mathcal{V}$, where $\|\cdot\|$ is a norm on $\mathbb{F}$ \emph{(}cf.~Section \ref{fdaFramework}\emph{)}.
\end{lemma}
Requiring that $D(x,P)\rightarrow \inf_{x\in\mathfrak{F}}D(x,P)$ as $\|x(v)\| \rightarrow \infty$ for Lebesgue almost every  $v\in \mathcal{V}$ is one natural analogue of Property 4.~of \citet{Zuo00, Liu90} and was suggested in \citet{Iwfos}, but we view Property P-3.~as more suitable in view of the arguments already set forth in this discussion. The second implication of P-3.~is Lemma \ref{lemmaDeepestPointUnique}. 

\begin{lemma}\label{lemmaDeepestPointUnique}
Let $D(\cdot,\cdot):\mathfrak{F}\times \mathcal{P}\rightarrow\mathbb{R}$ satisfy Property P-3 and let $z$ be as in P-3. Then $z=\argmax_{x\in\mathfrak{F}}D(x,P)$.
\end{lemma}

The direct analogue of Property 3 of \citet{Zuo00, Liu90} is to relax the strict inequality in Property P-3.
The strict inequality in P-3.~yields fewer ties in depth computations which enables us to better differentiate amongst the different elements of $\mathfrak{F}$. Moreover, strict inequality in P-3 automatically implies non-degeneracy of functional depth because it prevents all the points in $\mathfrak{F}$ having the same depth. Degenerate behaviour of several depth constructions is observed in \citet{IndioD}. They show that, inter alia, the band depth and half region depth constructions result in zero depth of every function in $\mathfrak{F}$ with probability one for common distributions such as continuous Gaussian processes.

\subsubsection{Discussion of P-4.~Upper semi-continuity in $x.$}\label{discussP4} In $\mathbb{R}$, there is a clear correspondence between the definition of depth and the cumulative distribution function $F(x)=P(X\leq x)$. The two natural ways of defining the depth at a point $x\in\mathbb{R}$ are $D(x,P)=P(X\leq x)\cdot P(X\geq x)$ and $D(x,P)=\min\{P(X\leq x),P(X\geq x)\}$, thus, from the c\`adl\`ag property of the cumulative distribution function, it is clear that, in $\mathbb{R}$, the depth is upper semicontinuous in the distance from the deepest point. The point we wish to make here is that, in order for the depth to reveal the features of the underlying distribution, it should, as a minimal requirement, satisfy the same properties as a cumulative distribution function, namely being non-decreasing (P-3.) and upper-semicontinuous (P-4.).

\subsubsection{Discussion of P-5.~Receptivity to convex hull width across the domain.}\label{discussP5} Many functional data sets encountered in practice contain functional data points that exhibit little variability over a particular subset of the domain $L\subset \mathcal{V}$, and significantly overlap with one another on $L$. The phenomenon described arises, inter alia, in functional microarray datasets \citep{CabreraLibro} and in chemometric datasets \citep[see e.g.~the yarn dataset in the R package pls][]{yarn}. Although the instinct is to draw parallels with the notion of heteroskedasticity in linear regression, this is in fact an entirely different phenomenon as it is usually still appropriate to view functional data as i.i.d.~copies of a random function $X$; $X$ simply possesses a variance function that is close to zero over $L$ and a correlation function close to one over $L\times L$. P-5.~obligates the depth to take heed of the values of $x\in\mathfrak{C}(\mathfrak{F},P)$ over $\mathcal{V}\backslash L$ to a greater extent than over $L$. Heuristically, the order of the curves does not matter much over $L$. Property P-5.~is particularly important in view of the discussion of P-6.~because, over $L$, small measurement error can conveibably lead to reconstructed functions that overlap in a drastically different way to the same functions observed without measurement error. A simple solution available for integrated depth constructions is to integrate the pointwise depths over a weight function depending on the convex hull of the data. This solution, proposed in \citet{HubertC}, effectively reduces the influence of regions over which all functions nearly coincide.

\subsubsection{Discussion of P-6.~Continuity in $P$.}\label{discussP6} Examples of $d_{\mathcal{P}}(\cdot,\cdot)$ are the Prohorov and bounded Lipschitz metrics, which both metricise the topology of weak convergence in the sense that 
$d_{\mathcal{P}}(Q,P)\rightarrow 0$ $P$-almost surely is equivalent to $Q\rightarrow P$ $P$-almost surely \citep[e.g.][Theorem 11.3.3]{DudleyRAP}.

\emph{Almost sure convergence of empirical depth to population depth}. The importance of Property P-6.~is evident when replacing $Q$ with $P_{n}$. In this case, fulfilment of P-6.~implies that the depth based on the empirical distribution converges almost surely to its population counterpart, i.e.~the estimator $D(\cdot,P_{n})\rightarrow D(\cdot,P)$ $P$-almost surely. This is particularly important when the depth is to be used for statistical inference. In this case, the objective is to gain understanding of population truths based on a random sample from that population. By contrast, in data analysis problems, the statistician typically has access to the whole population. Functional data analysis is however slightly different in view of the inherent partial observability of functional data.

\emph{Partial observability of functional data}. A second fundamental observation pertaining to P-6.~is that it tacitly addresses the inherent partial observability problem of functional data analysis. The latter gives rise to the delicate challenge of $P_{n}$ being inaccessible in its entirety. More specifically, whilst $P_{n}$ is a collection of weighted point masses at $X_{1},\ldots,X_{n}$, each valued in $\mathfrak{F}$, the practitioner only has access to $P_{n}'$, a collection of weighted point masses on $\{X_{i}(\mathcal{V}'_{i}): i=1,\ldots,n\}$, where $\mathcal{V}'_{i}\subset \mathcal{V}$ is a finite set that may be different for every $i\in\{1,\ldots,n\}$. 
The issue of partial observability of functional data is usually addressed through a preliminary interpolation or smoothing step to obtain an approximate reconstruction of the functional data object. Let $\widetilde{X}_{1},\ldots,\widetilde{X}_{n}$ be a sample of reconstructed functional data objects obtained from the random sample $\{X_{i}(\mathcal{V}'_{i}): i=1,\ldots,n\}$ of partially observed functional data objects or even from $\{X^{*}_{i}(\mathcal{V}'_{i}): i=1,\ldots,n\}$, where $X^{*}_{i}(\mathcal{V}'_{i})=\{X_{i}(v)+\epsilon_{i}, v\in\mathcal{V}'_{i}\}$ with $\{\epsilon_{i}:i=1,\ldots,n\}$ independent mean zero noise variables. Let $\widetilde{P}_{n}$ be the empirical probability measure over $\widetilde{X}_{1},\ldots,\widetilde{X}_{n}$. Then provided the reconstruction is such that $\widetilde{P}_{n}\rightarrow P$ $P$-almost surely, then Property P-6.~delivers the desired convergence of the functional depth. 

\emph{Qualitative robustness}. Importantly, fulfilment of P-6.~produces an embodiment of the empirical depth with the quintessential feature of robustness (cf.~Theorem \ref{thmHampel} below). The following definition of qualitative robustness is a restatement of Definition (A) in \citet{HampelAoS} in the more specific terms of the empirical depth. Here, we subscript the empirical depth by $P$ and $Q$ to emphasise that $P_{n}$ and $Q_{n}$ are random draws from $P$ and $Q$ respectively. With this notation, $\mathcal{L}(D_{P}(\cdot,P_n))$ is the probability measure on $\mathbb{R}$ induced by the mapping $D_{P}(\cdot,P_n)$. The theorem, and definition, are stated in terms of $P_n$ but it applies analogously when $P_{n}$ is replaced by $\widetilde{P}_{n}.$

\begin{definition}\label{definitionRobustness}\emph{[Qualitative robustness].} Let $P_{n}$ and $Q_{n}$ be the empirical measures corresponding to the $n$ random draws from $P$ and $Q$ respectively. For any $x\in\mathfrak{F}$, $D_{P}(x,P_n)$ is \emph{robust} at $P\in\mathcal{P}$ if and only if for all $\epsilon>0$ there exists a $\delta>0$ such that for any $Q\in\mathcal{P}$ satisfying 
$d_{\mathcal{P}}(Q,P)<\delta$, 
$d_{\mathcal{P}}(\mathcal{L}(D_{P}(x,P_n),\mathcal{L}(D_{Q}(x,Q_n)))<\epsilon$ for all $n.$
\end{definition}

\begin{theorem}\label{thmHampel}\emph{[An application of \citet{HampelBook}, Section 2.2, Theorem 2]}
If $D_{P}$ satisfies Property P-6, then $D_{P}(\cdot,P_n)$ is robust at $P$ for any $P\in\mathcal{P}$.
\end{theorem}

Qualitative robustness of the empirical depth is a desirable feature as it ensures that conclusions are not inordinately affected by outliers. 

\subsection{Implications for applications}\label{sectionApplications}

%\textcolor{red}{Certain properties amongst P-1.~to P-6.~are especially important for different types of applications, whilst others }

In this section we emphasise the roles played by P-1.~to P-6.~for different kinds of application.

Regarding P-1, in many applications, one would like the conclusions of statistical analysis or inference to be invariant to changes in the units of measurement. Nevertheless, for applications in which the ranking amongst the functions is the object of interest rather than the precise value of the depth, a weaker requirement may be sought: invariance of the ordering rather invariance of the depth values. This requirement would be suitable for constructing trimmed sample statistics by discarding the most extreme order statistics. There are applications in which the value of the depths themselves are of interest, and thus invariance in the precise sense of P-1.~is important. For instance, in certain model systems, systemic stability is related to diversity of a population and distance of the population centre of symmetry from a point, $p$, that is independent of the population. This situation arises in the model of the financial system considered by \citet{BealePNAS}, where $P=P_{n}$, i.e.~the whole population is available. One may construct a measure, $R$, of systemic risk from $d(z,p)$ and $\sum_{i=1}^{n}D(X_{i},P)$, where $z=\argmax_{x\in\mathfrak{F}} D(x,P)$, $P=P_{n}$ and $X_{i}$ is the relevant functional observation on individual $i$. The systemic risk contribution of individual $i$ is then $R_{i}=D(X_{i},P)/R$. This hints at the possibility of regulatory mechanisms designed to incentivise high systemic risk individuals towards a more systemically stable configuration in $\mathfrak{F}$ space. P-3.~ensures the diversity information is captured in the prescribed depths whilst P-3.~and P-4.~together ensure that the depth is not simply a ranking but captures the relative proximities of each individual to the centre of symmetry.

The centre-outward ordering induced by P-2.~and P-3.~and the information on relative proximities induced by P-3.~and P-4.~are qualities that enhance the ability of functional DD classifiers \citep{LiuCA2012} to differentiate between samples drawn from two different distributions. Moreover, the centre-outward ordering guaranteed by P-3.~provides the necessary and sufficient conditions for defining nearest neighbours \citep{Paindaveine2014}. Depth-based nearest neighbours have been effectively exploited (in the same reference) to define new classifiers, but they also offer prospects for nearest neighbour based nonparametric regression \citep[e.g][]{DevroyeGyorfiBook}.

Property P-5.~is also important for functional classification. If curves are from two different populations, both possessing covariance function close to zero over a subset $L\subset \mathcal{V}$, a functional DD plot classifier based on a depth violating P-5.,~ceteris paribus, has less power to discriminate between the two samples than one based on a depth satisfying P-5. Since classification is an inference (supervised learning) problem, P-6.~is important for ensuring that the sample depths of each $x\in \mathcal{F}$ converge to the corresponding population depths as $n\rightarrow \infty$. This assumption underpins the success of the DD classifier. 

Regardless of the precise nature of the application, P-6 is important for all of them, with its precise role depending on whether the application concerns inference or data analysis. For inference problems, the requirement is that the empirical depth converges to the population depth. Moreover, we require that the empirical depth based on the discretised functional data converges to the population depth. For data analysis problems, the aim is for the empirical depth based on the discretised functional data to converge to the empirical depth.

\section{A comparative study of existing functional depth proposals}\label{SectionCompareDepths}

In this section, we explore several popular constructions that have been proposed as functional depths in the literature. As we will see in due course, there is no single construction that satisfies all six properties in our definition of functional depth, which emphasises the necessity for further work in the area. Only functional depth constructions that have been proposed at the population level rather than simply at the sample level are explored, which rules out the construction based on distances that appears in \cite{Iwfos} and the one based on tilting that appears in \cite{GentonHall}.

\subsection{Existing functional depth constructions}\label{existingConstructions}

In each of the depth constructions outlined below, $X$ is a functional random variable defined on the probability space $(\mathfrak{F},\mathcal{A},P)$ (cf.~Section \ref{fdaFramework}) and, where relevant, expectation $\mathbb{E}$ is taken is with respect to $P$ unless otherwise stated. Sample analogues are obtained by replacing $P$ by $P_{n}$ for the idealised case and by $\widetilde{P}_{n}$ for the practically relevant case in which functional data objects are only observed at a finite set of evaluation points (cf.~Section \ref{discussP6}). For completeness, the sample versions of each depth construction in the idealised case are included after their population counterparts. The non-idealised sample versions, $D(\cdot,\widetilde{P}_{n}),$  are obtained by replacing $\{X_{1},\ldots,X_{n}\}$ by $\{\widetilde{X}_{1},\ldots,\widetilde{X}_{n}\},$ in $D(\cdot,P_{n})$. The constructions below need not uniquely prescribe a choice of metric $d$, however in most cases, there is a natural choice of $d$ with which to assess the fulfilment of Properties P-1.~to P-6.~in Definition \ref{definitionFDepth}. In each construction, $(\mathfrak{F},d)$ is as stated, $\mathcal{A}$ is the Borel sigma algebra (also the $d$-ball $\sigma$-algebra: cf.~Section \ref{fdaFramework}), and $P$ is a probability measure on the Borel sets of $\mathcal{A}$.

\subsubsection{The $h$-depth.} Let $(\mathfrak{F},d)=(\mathbb{H},\|\cdot\|_{\mathbb{L}_{2}})$. The $h$-depth \citep{Cuevas07} at $x\in \mathbb{H}$ w.r.t.~$P$ is defined as 
\begin{equation}\label{eqHDepth}
D_{h}(x,P):=\mathbb{E}K_{h}(\|x-X\|_{\mathbb{L}_{2}})
\end{equation}
where, for fixed $h>0$, $K_h(\cdot)=(1/h)K(\cdot/h),$ with $K(\cdot)$ the Gaussian kernel. The sample analogue of \eqref{eqHDepth} is $D_{h}(x,P_{n}):=\frac{1}{n}\sum_{i=1}^{n}K_{h}(\|x-X_{i}\|_{\mathbb{L}_{2}})$. The $h$-depth is the only example we consider that can be described as local (c.f.~Section \ref{discussP2}); a feature that is dispelled when the parameter $h$ is sufficiently large.

\subsubsection{The random Tukey depth.} Let $(\mathfrak{F},d)=(\mathbb{H},\|\cdot\|_{\mathbb{L}_{2}})$. Defining $\mathfrak{U}:=\{u_{1}, \ldots, u_{k}\}$, where $u_{j}$ $j=1,\ldots, k$ are realisations of $U_{j}$ $j=1,\ldots, k$, each drawn independently from a nondegenerate probability measure $\mu$ on $\mathbb{H}$, the random Tukey depth \citep{randomTukey} at $x\in\mathbb{H}$ w.r.t.~$P$ is
\begin{equation}\label{rTukey}
D_{RT}(x,P)=D_{\mathfrak{U}}(x,P):=\min_{u\in\mathfrak{U}}D_{1}(\langle u, x \rangle,P_{u})
\end{equation}
where, for any probability measure $Q$ on the Borel sets of $\RR,$ $D_{1}(t,Q)=\min\{Q(-\infty,t],Q[t,-\infty)\}$, $P_{u}$ is the marginal of $P$ on $\{\langle u,x\rangle:x\in\mathbb{H}\}$. $\mu$ is taken as a non-degenerate stationary Gaussian measure on $\mathbb{H}$.  For a discussion of the choice of $k$ see \citet{randomTukey}. The sample analogue of \eqref{rTukey} is simply obtained by replacing $P$ with $P_{n}$.

\subsubsection{The band depth.} Let $(\mathfrak{F},d)=(\mathcal{C},\|\cdot\|_{\infty})$ and let $\mathcal{V}\subset \mathbb{R}$. For $j\geq 2$, introduce the random $j$-simplex in $\mathfrak{F}$, $\mathcal{S}_{j}(P)=\{y\in\mathfrak{F}:y(v)=\alpha_{1}X_{1}(v)+\cdots +\alpha_{j}X_{j}(v): (\alpha_{k})_{k=1}^{j}\in \Delta^{j} \; \forall v\in\mathcal{V}, (X_{k})_{k=1}^{j}\sim P\}$, where $\Delta^{j}\subset \RR^{j-1}$ is the unit $j$-simplex. The band depth \citep{Romo09} at $x\in\mathfrak{F}$ is defined as
\begin{equation}\label{bandDepth}
D_{J}(x,P)=\sum_{j=2}^{J}P_{\mathcal{S}_j}\bigl(x\in \mathcal{S}_{j}(P)\bigr),
\end{equation}
where $P_{\mathcal{S}_j}$ is the probability measure over the random simplices constructed from the random $j$-tuple $X_{1},\ldots, X_{j}$.

When $P$ is replaced by $P_{n}$, there are $n$ choose $j$ distinct sets in the set of all random $j$-simplices on $\mathfrak{F}$ giving rise to the sample analogue of equation \eqref{bandDepth}, $D_J(x,P_n)=\sum_{j=2}^J {n\choose j}^{-1}\sum_{1\leq i_1<\ldots<i_j\leq n}\ind\bigl\{x\in B_{i_j}\bigr\}$, where 
$B_{i_{j}}:=\{y\in\mathfrak{F}:y(v)=\alpha_{1}X_{i_{1}}(v)+\cdots +\alpha_{j}X_{i_{j}}(v): (\alpha_{k})_{k=1}^{j}\in \Delta^{j} \; \forall v\in\mathcal{V}\}$ and $\{(i_{1},\ldots,i_{j}):i=1,\ldots,n\}$ defines the set of all possible $j$-tuples from $X_{1},\ldots,X_{n}$.

\subsubsection{The modified band depth.} Let $(\mathfrak{F},d)=(\mathcal{C},\|\cdot\|_{\infty})$ and let $\mathcal{V}\subset \mathbb{R}$. For $j\geq 2$, define a random $j$-simplex in $\RR$ to be of the form $\mathcal{S}_{j}(v,P)=\{y(v)\in\RR:y(v)=\alpha_{1}X_{1}(v)+\cdots +\alpha_{j}X_{j}(v): (\alpha_{k})_{k=1}^{j}\in \Delta^{j}, (X_{k})_{k=1}^{j}\sim P\}$, where $\Delta^{j}\subset \RR^{j-1}$ is the unit $j$-simplex. The modified band depth \citep{Romo09} at $x\in\mathfrak{F}$ is
\begin{equation}\label{modifiedBandDepth}
D_{MJ}(x,P)=\sum_{j=2}^{J}\mathbb{E}\bigl[\lambda \bigl\{v\in\mathcal{V}: x\in \mathcal{S}_{j}(v,P)\bigr\}\bigr]/\lambda(\mathcal{V}),
\end{equation}
where expectation is with respect to the measure $P_{\mathcal{S}_j}$, as defined above in the definition of the band depth. In Section \ref{sectionProofs} it will sometimes be convenient to refer to $S_{j}(v,P)=[L_{j}(v),U_{j}(v)]$, where $L_{j}(v):=\min_{y\in \mathfrak{X}_{j}}y(v)\}$ and $U_{j}(v):=\max_{y\in \mathfrak{X}_{j}}y(v)$, where $\mathfrak{X}_{j}=(X_{1},\ldots,X_{j})$ and $X_{1},\ldots,X_{j}\sim P$.

When $P$ is replaced by $P_{n}$, there are $n$ choose $j$ distinct sets in the set of all random $j$-simplices on $\mathfrak{F}$ giving rise to the sample analogue of equation \eqref{modifiedBandDepth}, 
\[
D_{MJ}(x,P_{n}):= \sum_{j=2}^J{n\choose j}^{-1}\sum_{1\leq i_1<\ldots<i_j\leq n}\lambda\bigl\{v\in\mathcal{V}: x(v) \in B_{i,j}(v)\bigr\}/\lambda(\mathcal{V}),
\]
where $B_{i_{j}}(v):=\{y(v)\in\RR:y(v)=\alpha_{1}X_{i_{1}}(v)+\cdots +\alpha_{j}X_{i_{j}}(v): (\alpha_{k})_{k=1}^{j}\in \Delta^{j}\}$ and $\{(i_{1},\ldots,i_{j}):i=1,\ldots,n\}$ defines the set of all possible $j$-tuples from $X_{1},\ldots,X_{n}$.

\subsubsection{The half region depth.} In the same setting as for the band depth, the half region depth \citep{LopezHalfRegion} w.r.t.~$P$ at $x\in \mathfrak{F}$ is
\begin{equation}\label{eqHR}
D_{HR}(x,P):=\min\{P(X\in H_{x}),P(X \in E_{x})\}.
\end{equation}
where $H_{x}$ is the hypograph of $x$, i.e.~$H_{x}:=\{y\in\mathfrak{F}:y(v)\leq x(v) \; v\in\mathcal{V}\}$, and $E_{x}$ is the epigraph of $x$, i.e.~$E_{x}:=\{y\in\mathfrak{F}:y(v)\geq x(v) \; v\in\mathcal{V}\}$. Thus the halfspace depth is the minimum between the proportion of curves in the epigraph and hypograph of $x$. The sample analogue of \eqref{eqHR} is obtained by replacing $P(X\in H_{x})$ in \eqref{eqHR} by $\frac{1}{n}\sum_{i=1}^{n}\ind\bigl\{X_{i}\in H_{x}\bigr\}$ and analogously for $P(X \in E_{x})$.

\subsubsection{The modified half region depth.} In the same setting as for the band depth, the half region depth \citep{LopezHalfRegion} w.r.t.~$P$ at $x\in \mathfrak{F}$ is
\[
D_{MHR}(x,P)=\min\Bigl\{\mathbb{E}[\lambda\{v\in \mathcal{V}: X(v)\leq x(v)\}],\mathbb{E}[\lambda\{v\in \mathcal{V}: X(v)\geq x(v)\}]\Bigr\}/\lambda(\mathcal{V}),
\]
with sample analogue
\[
D_{MHR}(x,P_n)=\min\Bigl\{\frac{1}{n}\sum_{i=1}^n\lambda\{v\in \mathcal{V}: X_i(v)\leq x(v)\}, \frac{1}{n}\sum_{i=1}^n\lambda\{v\in \mathcal{V}:X_i(v)\geq x(v)\}\Bigr\}/\lambda(\mathcal{V}).
\]
In Table \ref{tableDepthDefinitions}, we summarise the depth constructions presented in detail above. 

\begin{table}
\caption{Summary of existing depth constructions.}
\begin{center}
\begin{tabular}{|c|c|c|c|}
\hline
Depth & ($\mathfrak{F},d$) & $\mathcal{V}$ & Construction \\
\hline 
$D_{h}(x,P)$      & $\mathbb{H}$, $\|\cdot\|_{\mathbb{L}_{2}}$ & $\mathcal{V}\subset \mathbb{R}^{d}$, $d\geq 1$ & $\mathbb{E}K_{h}(\|x-X\|_{\mathbb{L}_{2}})$\\
$D_{RT}(x,P)$     & $\mathbb{H}$, $\|\cdot\|_{\mathbb{L}_{2}}$ & $\mathcal{V}\subset \mathbb{R}^{d}$, $d\geq 1$ & $\min\bigl\{D_{1}(\langle u, x \rangle,P_{u}): u \in \mathfrak{U}\bigr\}$ \\
$D_{J}(x,P)$      & $\mathcal{C}$, $\|\cdot\|_{\infty}$          & $\mathcal{V}\subset \mathbb{R}$ & $\sum_{j=2}^{J}P_{\mathcal{S}_j}(x\in \mathcal{S}_{j}(P))$\\
$D_{MJ}(x,P)$     & $\mathcal{C}$, $\|\cdot\|_{\infty}$          & $\mathcal{V}\subset \mathbb{R}$ & $\sum_{j=2}^{J}\mathbb{E}[\lambda\{v\in\mathcal{V}: x(v)\in \mathcal{S}_{j}(v,P)\}]/\lambda(\mathcal{V})$\\
$D_{HR}(x,P)$     & $\mathcal{C}$, $\|\cdot\|_{\infty}$          & $\mathcal{V}\subset \mathbb{R}$ 								&$\min\{P(X\in H_{x}),P(X\in E_{x})\}$\\
$D_{MHR}(x,P)$    & $\mathcal{C}$, $\|\cdot\|_{\infty}$          & $\mathcal{V}\subset \mathbb{R}$ 								& $\min\{\mathbb{E}[\lambda\{v\in \mathcal{V}: X(v)\leq x(v)\}]/\lambda(\mathcal{V}),$\\
                  &                                              &               &  \qquad  \;    $\mathbb{E}[\lambda\{v\in \mathcal{V}: X(v)\geq x(v)\}]/\lambda(\mathcal{V})\}$\\
\hline
\end{tabular}
\end{center}
\end{table}\label{tableDepthDefinitions}

\subsubsection{Other existing functional depth proposals}

In addition to the six functional depth proposals exposed above, there are several other constructions that have appeared in the literature. The integrated depth is proposed in \citet{Fraiman} as the first depth for functional data. It is defined by integrating over the continuum of one dimensional pointwise depths at each point $x(v)$, $v\in\mathcal{V}$. As noted in \citet{HubertC}, the integrated depth is related to the modified band depth of \citet{Romo09}. More specifically, the modified band depth with $J=2$, the recommended value in \citet{Romo09}, coincides with the integrated depth when computed w.r.t.~a probability distribution with absolutely continuous marginals. This correspondence is due to the use of the simplicial depth for the one dimensional pointwise depth, as initially proposed in \citet{Fraiman}. Other one dimensional pointwise depths are equally valid, but do not give rise to this same link with the modified band depth. The multivariate functional halfspace depth of \citet{HubertC} generalises the integrated depth, allowing multivariate functions through the use of the multidimensional pointwise Tukey depth, and through the inclusion of a weight function to downweight the influence of the pointwise depth values over regions where the convex hull width is small. Another approach to generalise the integrated depth to multivariate functions that was proposed in the recent literature is in \citet{Hlubinka}. Other functional depth proposals include the integrated dual depth of \citet{Cuevas09}, proposed as the population analogue of the random projection depth \citep{Cuevas07}. There, the double random projection depth was also proposed as the first example of depth suitable for multivariate functional data. Additionally, \citet{ChakrabortyAoS2014} and \citet{IndioD} study from a functional depth perspective the spatial depth of \citet{Chaudhuri1996}, \citet{VardiZhang2000} and \citet{Serfling2002} to the functional setting, which has proved to be a useful construction. For a generalisation of some of these depths, see \citet{Mosler2013}.

\subsection{A property-wise analysis of existing functional depths}

In the theoretical results that follow, $D_{h}$, $D_{RT}$, $D_{J}$, $D_{MJ}$, $D_{HR}$, $D_{MHR}$ and their respective $(\mathfrak{F},d)$ are as in Table \ref{tableDepthDefinitions}, and $\mathcal{D}:=\{D_{h},D_{RT},D_{J},D_{MJ},D_{HR},D_{MHR}\}$. The conclusions of the following theorems are summarized in Table \ref{tableprop}. We comment here on reasons for which the different examples of depth satisfy, or fail to satisfy, the corresponding properties. For a deeper insight, see the proofs in Section \ref{sectionProofs}.

\begin{table}
\caption{Adherence of existing depth constructions to depth defining properties.}
\begin{center}
\begin{tabular}{|r|cccccc|}
\hline
             &P-1        &	P-2G    & P-3 & P-4 & 	P-5  & P-6  \\
\hline
$D_{h}$      &\tickNo	   &\tickYes  &\tickYes &\tickYes  &\tickYes& \tickYes \\
$D_{RT}$     &\tickYes	 &\tickYes  &\tickNo  &\tickYes  &\tickNo & \tickYes \\
$D_{J}$      &\tickYes	 &\tickYes  &\tickNo  &\tickYes  &\tickNo & \tickYes  \\
$D_{MJ}$     &\tickYes	 &\tickYes  &\tickNo  &\tickYes  &\tickNo & \tickYes  \\
$D_{HR}$     &\tickYes	 &\tickNo   &\tickNo  &\tickYes  &\tickNo & \tickYes  \\
$D_{MHR}$    &\tickYes	 &\tickYes  &\tickNo  &\tickYes  &\tickNo & \tickYes \\
\hline
\end{tabular}
\end{center}
\end{table}\label{tableprop}

\begin{theorem}\label{thmP-1}\emph{[Property P-1.~Distance invariance]}. All elements of $\mathcal{D}$ satisfy Property P-1 with the exception of $D_{h}$.
\end{theorem}
The part of the proof of Theorem \ref{thmP-1} concerning the $h$ depth assumes that the same $h$ is used in $D_{h}(x,P_{X})$ and $D_{h}(f(x),P_{f(X)})$, but the conclusion remains valid if we allow for $h$ to depend on $f$. To see it, simply observe that $\frac{1}{h}\exp\{-\|x-X\|^2/2h^{2}\} \neq \frac{1}{h_{f}}\exp\{-a\|x-X\|^2/2h_{f}^{2}\}$ for any $h>0, h_{f}>0.$ 

Recall from our discussion of P-2.~that, since there is no unique measure of centre of symmetry, $\theta$, in general, it is more meaningful to consider the behaviour of $D$ for a particular case of $P$ in which all standard notions of centre of symmetry coincide at $\theta$. We thus consider here adherence to P-2G. 

\begin{theorem}\label{thmP-2}\emph{[Property P-2G.~Maximality at Gaussian process mean]}. 
With the exception of $D_{HR}$, all elements of $\mathcal{D}$ satisfy Property  P-2G, where $J\geq 3$ in $D_{J}.$
\end{theorem}

The intuitive explanation for $D_{HR}$ failing to satisfy P2-G is that the expected number of upcrossings of a mean zero Gaussian process above a level $a$ is strictly decreasing in $|a|$. Hence the probability that a Gaussian process is either entirely above or entirely below $a$ is strictly increasing in $|a|$. The modified version of $D_{HR}$ does not suffer this drawback as it takes account of the duration of excursions above $|a|$.

For sufficiently small $h$, the $h$-depth becomes a local depth rather than a global depth, and hence, as alluded to in the discussion in Section \ref{discussP3}, one would not expect a centre outward ordering from a unique centre of symmetry, but rather an outward ordering from points of high local depth. As such, verification of P-3 is only achievable when $h$ is sufficiently large for $D_{h}$ to constitute a global depth. We implicitly impose this assumption in Lemma \ref{P-3h} below by imposing that the deepest element (as measured by $D_{h}$) exists and coincides with the mean.

\begin{lemma}\label{P-3h}
Provided that $\mathbb{E}X$ exists and $D_h(\mathbb{E}X,P)=\sup_{x\in\mathfrak{F}}D_h(x,P)$, $D=D_{h}$ satisfies P-3.
\end{lemma}

Lemma \ref{P-3h} works for any type of distribution, including both continuous and discrete. However, the counterexamples in the proof of Theorem \ref{thmP-3} demonstrate that non-continuous distributions preclude adherence to P-3.~for elements of $\mathcal{D}\setminus\{D_h\}$. The constructions of these depths are based more directly on terms of the form $P(B_{x})$ for $B_{x}$ a Borel set that depends on $x\in\mathfrak{F}$. For non-continuous distributions and the constructions we consider, there exist $x,y\in \mathfrak{F}$ with $x\neq y$ that yield $P(B_{x})=P(B_{y})$, resulting in the assignment of equal depths to $x$ and $y$.

\begin{theorem}\label{thmP-3}\emph{[Property P-3.~Strictly decreasing w.r.t.~the deepest point]}.
The elements of $\mathcal{D}\setminus\{D_h\}$ do not satisfy Property  P-3.
\end{theorem}

Lemma \ref{lemmaHContinuous}, as well as being of independent interest, is used in the proof of Theorem \ref{thmP-4}.

\begin{lemma}\label{lemmaHContinuous}
For any $P\in\mathcal{P}$, $D_{h}(x,P)$ is continuous in $x$.
\end{lemma}

\begin{theorem}\label{thmP-4}\emph{[Property P-4.~Upper semi-continuity in $x$]}.
All elements of $\mathcal{D}$ satisfy Property  P-4. 
\end{theorem}

Upper semi-continuity of the elements of $\mathcal{D}$ arises naturally because all depth constructions preserve the upper semi-continuity of the distribution function induced by $P$. A stricter requirement of continuity would, in most cases, rule out the possibility of $P$ with finite support.

\begin{theorem}\label{thmP-5}\emph{[Property P-5.~Receptivity to convex hull width across the domain]}.
Provided that $\mathbb{E}X$ exists, $D=D_{h}$ satisfies P-5. The elements of $\mathcal{D}\setminus\{D_h\}$ do not satisfy Property  P-5.
\end{theorem}

The intuition behind the non-adherence of the elements of $\mathcal{D}\backslash D_{h}$ to P-5.~is that their constructions all result in an assignment of rank, neglecting the relative distances (as measured in some suitable metric, $d,$ with respect to $P$) between elements of $\mathfrak{F}$. By contrast, the $h$-depth is essentially a weighted $\LL_{2}(\mathcal{V},\lambda)$, where the weights depend on $P$. As such, it is able to appropriately exploit the information contained in $P$ such that the influence of variations in $X$ over $L_{\delta}$ is commensurate with $\delta$.

\begin{theorem}\label{thmP-6}\emph{[Property P-6.~Continuity in $P$]}.
All elements of $\mathcal{D}\backslash \{D_{J},D_{RT}\}$ satisfy Property  P-6. $D_{RT}$ satisfies P-6.~when the limiting distribution is continuous or the sequence of distributions is the sequence of empirical distributions. $D_{J}$ satisfies P-6.~when $\mathfrak{F}$ is restricted to be the space of equicontinuous functions on $\mathcal{V}\subset \RR$.
\end{theorem}

All elements of $\mathcal{D}\backslash\{D_{J},D_{MJ}\}$ are either constructed from sets of the form $P(B_{x})$ for $B_{x}$ a Borel set that depends on $x\in\mathfrak{F}$, or as an integral of a bounded Lipschitz function with respect to $P$, which yields adherence to P-6.~by the well known Portmanteau theorem for weak convergence (cf.~Section \ref{sectionProofs} for details). The construction of $D_{J}$ and $D_{MJ}$ results in a stochastic process whose behaviour is governed by $P$. As is shown in Section \ref{sectionProofs} convergence of $Q$ to $P$ guarantees weak convergence of the respective stochastic processes which in turn results in pointwise $P$-a.s.~convergence of depths.

Amongst the six constructions we consider, the $h$-depth satisfies 5 of the 6 properties we seek. This should not be interpreted as a recommendation to favour the $h$-depth. As discussed in section \ref{sectionApplications} each property has different implications for different application areas and a depth construction should thus be chosen with the application in mind. As the $h$-depth fails to satisfy P-1. A proposal is to substitute the proposed kernel. As a simple illustration, if the kernel function resulted in $D_{h}(x,P):=\frac{1}{\sqrt{2\pi}}\exp\{-\|x-X\|^2/2h^{2}\},$ property P-1 would be satisfied when allowing $h$ to depend on $f,$ where $f$ is defined in Definition \ref{definitionFDepth}.

\section{Proofs}\label{sectionProofs}

\begin{proof}[\scshape{Proof of Lemma \ref{lemmaLpSpaces}}.]
For any $x,z\in\mathfrak{F},$ $(d(x,z))^p\leq (\sup_{v\in\mathcal{V}}\|x(v)-z(v)\|)^p\lambda(\mathcal{V}).$ Fixing $z,$ as $\lambda(\mathcal{V})$ is finite, $d(x,z)\rightarrow\infty$ implies that $\sup_{v\in\mathcal{V}}\|x(v)\|\rightarrow\infty.$ Thus, $D(x,P)\rightarrow\inf_{x\in\mathfrak{F}}D(x,\mathfrak{F})$ as $\sup_{v\in\mathcal{V}}\|x(v)\|\rightarrow\infty$ and, a fortiori, as $\|x(v)\|\rightarrow\infty$ for Lebesgue almost every $v\in\mathcal{V}.$
\end{proof}

\begin{proof}[\scshape{Proof of Lemma \ref{lemmaDeepestPointUnique}}.]
Suppose for a contradiction that there exist $z_{1}, z_{2}\in\mathfrak{F}$ $z_{1}\neq z_{2}$ such that $D(z_{1},P)=D(z_{2},P)=\max_{x\in\mathfrak{F}}D(x,P)$. As $z_{1}\neq z_{2}$ implies $d(z_{1},z_{2})>0$, we may take in the statement of P-3.~$x=z_{1}$ and $z=z_2$, which yields by P-3.~$D(z_{1},P)<D(z_{2},P)$, a contradiction.
\end{proof}

\begin{proof}[\scshape{Proof of Theorem \ref{thmP-1} (Property P-1.)}.]

\emph{h-depth}. When $(\mathfrak{F},d)=(\mathbb{H},\|\cdot\|_{\LL_{2}})$, the set of functions that satisfy $d(f(x),f(y))=a_f \cdot d(x,y)$ for any $x,y\in \mathfrak{F}$ is given by 
\begin{equation}\label{classL2Functions}
\bigl\{f : f(x(v))=\sqrt{a(v)}x(v), \; a(v)=a_f>0 \;\;\forall v\in \mathcal{V}\bigr\}.
\end{equation}
Since $K_h(a\|x-X\|)\neq K_h(\|x-X\|)$ for all $a\neq 1$, there exist functions in the set \eqref{classL2Functions} for which $D_{h}(x,P_X) \neq D_{h}(f(x),P_{f(X)})$.
 \emph{Random Tukey depth}: Let $(\mathfrak{F},d)=(\mathbb{H},|\langle \cdot,\cdot\rangle|)$, then the set of functions that satisfy $d(f(x),f(y))=a_f \cdot d(x,y)$ for any $x,y\in \mathfrak{F}$ is given by equation \eqref{classL2Functions}. The result follows since $\{y:\langle u, x-y\rangle \geq 0\}=\{y:\langle u,\sqrt{a}y-\sqrt{a}x\rangle\geq 0\}$ for all $v\in\mathbb{H}$.
  
For $D_{J}$, $D_{MJ}$, $D_{HR}$ and $D_{MHR},$ let $(\mathfrak{F},d)=(\mathcal{C}(\mathcal{V}),\|\cdot\|_{\infty}),$ the set of functions satisfying $d(f(x),f(y))=a_f \cdot d(x,y)$ for any $x,y\in \mathfrak{F}$ is given by 
\[
\bigl\{f:f(x(v))=a(v)x(v)+b(v), \; |a(v)|=a_f>0 \;\; \forall v\in\mathcal{V}\bigr\}.
\]
Then, $D(x,P_X)= D(f(x),P_{f(X)})$ for those instances of depth listed above by  the following observations.  
\emph{Band depth}: the result is Theorem 3 of \citet{Romo09}. \emph{Modified band depth}: for $a_f>0,$ $\displaystyle x(v)\in[L_j(v),U_j(v)]$ if and only if $ a_f x(v)\in[a_f L_j(v),a_f U_j(v)]$. \emph{Half-region depth}: we have $P[X(v) \leq x(v), v \in \mathcal{V}]=P[a_f X(v) \leq a_f x(v), v \in \mathcal{V}]$ and analogously for $P[X(v) \geq x(t), v \in \mathcal{V}]$. \emph{Modified half-region depth}: we have $E[\lambda\{v \in \mathcal{V}: x(v) \leq X(v)\}]=E[\lambda\{v \in \mathcal{V}: a_f x(v) \leq a_f X(v)\}]$ and analogously for $E[\lambda\{v\in \mathcal{V} : x(v) \geq X(v)\}]$. 
\end{proof}

\begin{proof}[\scshape{Proof of Theorem \ref{thmP-2}} (Property P-2G.).]

\emph{$h$-depth}. Suppose for a contradiction that $z:=\argsup_{x\in\mathfrak{F}}D_{h}(x,P)$ is such that $D_h(z,P_X)>D_h(E[X],P_X)$. Since
\begin{equation}\label{eqSupInf}
\argsup_x D_h(x,P_X)=\argsup_x \EE\Bigl[\exp\Bigl\{-\frac{\|x-X\|^2}{2h^2}\Bigr\}\Bigr]=\arginf_x \EE[\|x-X\|^2],
\end{equation}
the previous supposition is equivalent to $\EE[\|z-X\|^2]<\EE[\|\mathbb{E}[X]-X\|^2]=\EE[\|X\|^2]$. After some algebra we obtain $\|x\|^2<2\int x(v)\EE[X(v)]dv=0$, a contradiction.

\emph{Random Tukey depth}. For any $u\in\mathfrak{F}=\mathbb{H}$, we have that $\langle u, \mathbb{E}X \rangle$ is the mean of $P_{u}$ because $\mathbb{E}[\langle X,u \rangle]=\mathbb{E}\int X(v) u(v) dv = \int \mathbb{E}X(v) u(v) dv$. Since, for $P$ a Gaussian process, the mean of $P_{u}$ coincides with the median of $P_{u}$, we have $D_{1}(\langle u,\mathbb{E}X\rangle,P_{u})=\frac{1}{2}$. Then, by the definition of random Tukey depth, $D_{RT}(\mathbb{E}X,P)=\min_{u\in\mathfrak{U}}\frac{1}{2}=\frac{1}{2}$, the maximum attainable value for the random Tukey depth, hence $D_{RT}(\mathbb{E}X,P)=\sup_{x\in\mathfrak{F}} D_{RT}(x,P)$.

\emph{Band depth and modified band depth}. By the definition of the band depth and the modified band depth
\[
\sup_{x\in\mathfrak{F}}D_{J}(x,P)\leq \sum_{j=2}^{J}\sup_{x\in\mathfrak{F}}P_{\mathcal{S}_j}(x\in\mathcal{S}_{j}(P))
\]
and
\[
\sup_{x\in\mathfrak{F}}D_{MJ}(x,P)\leq \sum_{j=2}^{J}\sup_{x\in\mathfrak{F}}\mathbb{E}[\lambda\{v\in\mathcal{V}: x(v)\in\mathcal{S}_{j}(v,P)\}/\lambda(\mathcal{V})].
\]
Since each of $X_{1},\ldots, X_{J}$ is a random draw from $P$, whose mean is $\theta=\mathbb{E}X$, and since $P_{\mathcal{S}_j}$ is a continuous distribution over simplices (because $P$ is continuous), the $x$ which maximises the probability of a random $j$-simplex enveloping it is clearly $x=\theta$, yielding $\sup_{x\in\mathfrak{F}}D_{J}(x,P)=D_{J}(\theta,P)$. Similarly, the $x$ for which the expectated duration spent in any simplex is largest is also $x=\theta$, yielding $\sup_{x\in\mathfrak{F}}D_{MJ}(x,P)=D_{MJ}(\theta,P)$.

\emph{Half region depth}. By \citet{Adler1981} Theorem 4.1.1, the expected number of upcrossings of a level $\bar{u}$ of a zero-mean, stationary, almost surely continuous random process on $\mathcal{V}$, is 
\begin{equation}\label{eqnAdler}
\EE[N_{\bar{u}}]= \sqrt{\frac{-R''(0)}{R(0)}}\frac{\lambda(\mathcal{V})}{2\pi}\exp\left\{-\frac{\bar{u}^{2}}{2R(0)}\right\}
\end{equation}
where $R(0)=\EE[|X(v)|^{2}]$ and $-R''(0)$ is the variance of $X(v)$, which is constant by stationarity of $X$. Equation \eqref{eqnAdler} is maximised at $\bar{u}=0$, hence for any $\bar{u}$ such that $0<|\bar{u}|<\infty$, $\min\left\{P\bigl(X(v)\leq \bar{u} \; \forall v\in\mathcal{V}\bigr), P\bigl(X(v)\geq \bar{u} \; \forall v\in\mathcal{V}\bigr)\right\} > \min\left\{P\bigl(X(v)\leq 0 \; \forall v\in\mathcal{V}\bigr), P\bigl(X(v)\geq 0 \; \forall v\in\mathcal{V}\bigr)\right\}$.

\emph{Modified half region depth}. Demonstrating that $D_{HR}(x,P)$ achieves its maximum value at the zero mean function of the Gaussian process $P$, entails a proof that the expected measure of the level zero excursion set is 1/2, where the level zero excursion set is defined as
\[
A_{0}:=A_{0}(X,\mathcal{V}):=\{v\in\mathcal{V}: X(v)\geq 0\}.
\]
By \citet{Rice1945}, from which equation \eqref{eqnAdler} also originally derived, the expected length of an excursion above zero is $\pi\sqrt{R(0)/[-R''(0)]}$. Recalling that $\mathcal{V}$ is a compact subset of $\RR$ and assuming an excursion starts at $\min\{v\in \mathcal{V}\}$, we thus have, using equation \eqref{eqnAdler},
\[
\EE[\lambda(A_{0})]=\frac{\lambda(\mathcal{V})}{2}\sqrt{\frac{-R''(0)}{R(0)}}\sqrt{\frac{R(0)}{-R''(0)}}=\frac{\lambda(\mathcal{V})}{2}.
\]
Hence $D_{MHR}(\EE X, P)=1/2$, which coincides with $\sup_{x\in\mathfrak{F}}D(x,P)$.
\end{proof}

\begin{proof}[\scshape{Proof of Lemma \ref{P-3h}} (Property P-3.~\emph{$h$-depth}).]

Observe that $D_h$ is translation invariant, i.e.~for any  $x,b\in\mathfrak{F}$ and $P_X\in\mathcal{P}$
\[
D_h(x,P_X)=\EE\Bigl[\frac{1} {h\sqrt{2\pi}}exp\Bigl\{-\frac{\|x-X\|^2} {2h^2}\Bigr\}\Bigr]=\EE\Bigl[\frac{1} {h\sqrt{2\pi}}exp\Bigl\{-\frac{\|(x-b)-(X-b)\|^2} {2h^2}\Bigr\}\Bigr].
\]
Thus set $\EE[X]=0$ without loss of generality.

Suppose for a contradiction $D_h(x,P)\geq D_h(y,P)$. Substituting $\|x-X\|^2=\|x\|^2+\|X\|^2-2\int x(v)X(v)dv$ in the expression for $D_{h}$ gives the inequality
\[
\exp\Bigl\{-\frac{\|x\|^2-\|y\|^2}{2h^2}\Bigr\}\geq \EE\Bigl[\exp\Bigl\{\frac{\int(y(v)-x(v))X(v)dv}{h^2}\Bigr\}\Bigr].
\] 
By the statement of P-3 and the fact that $\EE[X]=0$, we have $\|x\|>\|y\|$ and so
\begin{eqnarray}\label{left}
1>\exp\Bigl\{-\frac{\|x\|^2-\|y\|^2}{2h^2}\Bigr\}.
\end{eqnarray}
On the other hand, by Jensen's inequality 
\[
\EE\Bigl[\exp\Bigl\{\frac{\int(y(v)-x(v))X(v)dv}{h^2}\Bigr\}\Bigr]\geq \exp\Bigl\{\frac{\int(y(v)-x(v))\EE[X](v)dv}{h^2}\Bigr\}
\]
which is equal to $1$ because $\EE[X]=0$. This together with (\ref{left}) yields the contradiction
\[
1>\exp\{-\frac{\|x\|^2-\|y\|^2}{2h^2}\}\geq 1.
\]
\end{proof}

\begin{proof}[\scshape{Proof of Theorem \ref{thmP-3}} (Property P-3.).]

\emph{Random Tukey depth}. The proof is by counterexample. Let $P\in\mathcal{P}$ be a discrete distribution with support $\{x_{1},x_{2}\}$ with $x_{1}(v)=2$ for all $v\in\mathcal{V}$ and $x_{2}(v)=-1$ for all $v\in\mathcal{V}$. Let $u\in\mathbb{H}$ be an arbitrary realisation of the random variable $U$ whose distribution is $\mu$. The inner product with $u$ of any $y\in\mathcal{Y}:=\{y(v)=c \; \forall v\in\mathcal{V} \text{ with } c\in(-1,2)\}$ gives rise to $\langle u,y \rangle\in(\min\{\langle u,x_{1}\rangle,\langle u,x_{2}\rangle\},\max\{\langle u,x_{1}\rangle,\langle u,x_{2}\rangle\})$. It follows that $D_{RT}(y,P)=\max_{x\in\mathfrak{F}}D_{RT}(x,P)$ for any $y$ in the closure of $\mathcal{Y}$, which contradicts Lemma \ref{lemmaDeepestPointUnique}.

\emph{Band depth}. The proof is by counterexample. Take $P\in\mathcal{P}$ discrete with $P(\{x_{1}\})=P(\{x_{2}\})=1/2$, where $x_{1}(v)=-c$ for all $v\in\mathcal{V}$, $x_{2}(v)=c$ for all $v\in\mathcal{V}$. Then $P_{\mathcal{S}_j}$ $j=J=2$ is discrete with $P_{\mathcal{S}_j}(\mathcal{S}_{j,1})=P_{\mathcal{S}_j}(\mathcal{S}_{j,2})=1/4$ and $P_{\mathcal{S}_j}(\mathcal{S}_{j,3})=1/2$, where $\mathcal{S}_{j,1}=\{x_{1}\}$, $\mathcal{S}_{j,2}=\{x_{2}\}$ and $\mathcal{S}_{j,3}=\{[x_{1}(v),x_{2}(v)]: v\in\mathcal{V}\}$. Then $D_{J}(z,P)$ has two global maxima, at $z=x_{1}$ and at $z=x_{2}$, with $D_{J}(z,P)=3/4$. Without loss of generality, set $z=x_{1}$. For any $x,y\in\mathfrak{F}=\mathcal{C}(\mathcal{V})$ such that $\max\{d(y,z),d(y,x)\}<d(x,z)$ and $x_{2}(v)<x(v)<x_{1}(v)$, $x_{2}(v)<y(v)<x_{1}(v)$ for all $v\in\mathcal{V}$. Then $D_{J}(x,P)=D_{J}(y,P)=1/2$, violating P-3.

\emph{Modified band depth}. The proof uses the same counterexample as in the proof for the band depth. We have
\begin{eqnarray*}
D_{MJ}(z,P)&=&\lambda\{v\in\mathcal{V}: z(v)\in \mathcal{S}_{j,1}(v,P)\}P_{\mathcal{S}_{j}}(\mathcal{S}_{j,1})/\lambda(\mathcal{V}) \\
           & &  + \;\; \lambda\{v\in\mathcal{V}: z(v)\in \mathcal{S}_{j,2}(v,P)\}P_{\mathcal{S}_{j}}(\mathcal{S}_{j,2})/\lambda(\mathcal{V}) \\
					 & &  + \;\; \lambda\{v\in\mathcal{V}: z(v)\in \mathcal{S}_{j,3}(v,P)\}P_{\mathcal{S}_{j}}(\mathcal{S}_{j,3})/\lambda(\mathcal{V}), 
\end{eqnarray*}
and $D_{MJ}(z,P)$ is maximised at $z=x_{1}$ and $z=x_{2}$, giving $D_{MJ}(z,P)=3/4$. Without loss of generality, set $z=x_{1}$. For any $x,y\in\mathfrak{F}=\mathcal{C}(\mathcal{V})$ such that $\max\{d(y,z),d(y,x)\}<d(x,z)$ and $x_{2}(v)<x(v)<x_{1}(v)$, $x_{2}(v)<y(v)<x_{1}(v)$ for all $v\in\mathcal{V}$. Then $D_{MJ}(x,P)=D_{MJ}(y,P)=1/2$, violating P-3.

\emph{Half region depth}. Let $P$, $x$ and $y$ be as for the (modified) band depth. Then $D(z,P)=P(X(v)\geq z(v), v\in\mathcal{V})=P(X(v)\leq z(v), v\in\mathcal{V})$. But $P(X(v)\geq x(v), v\in\mathcal{V})=P(X(v)\geq z(v), v\in\mathcal{V})=P(X(v)\geq y(v), v\in\mathcal{V})$ hence $D_{HR}(x,P)=D_{HR}(y,P)=D_{HR}(z,P)$ despite the fact that $d(y,z)<d(x,z)$.

\emph{Modified half region depth}. Let $P$, $x$ and $y$ be as for the (modified) band depth. Then for any $\omega\in\Omega$, $\lambda\{v\in\mathcal{V}:X(\omega,v)\leq x(v)\}=\lambda\{v\in\mathcal{V}:X(\omega,v)\leq y(v)\}$ and likewise for the converse inequality. Hence $D_{MHR}(x,P)=D_{MHR}(y,P)$ despite the fact that $d(y,z)<d(x,z)$.
\end{proof}

\begin{proof}[\scshape{Proof of Lemma \ref{lemmaHContinuous}}.]
Write $\exp\bigl\{-\|x-X(\omega)\|/2h\bigr\}/\sqrt{2\pi h}=:F(x,\omega)$.
Then for $P$-almost every $\omega\in\Omega$, $F(\cdot,\omega)$ is continuous at $x$. Moreover, since $\exp\{-z\}$ is bounded on $z\in\mathbb{R}^{+}$, there exists a $P$-integrable function $g(\omega)$ such that $F(y,\omega)\leq g(\omega)$ for $P$-almost every $\omega\in\Omega$ and all $y$ in a neighbourhood of $x$. Since the above holds for all $x\in\mathfrak{F}$, it follows by Theorem 7.43 of \citet{SPBook} that $\mathbb{E}\bigl[\exp\bigl\{-\|\cdot-X(\omega)\|/2h\bigr\}/\sqrt{2\pi h}\bigr]$ is continuous at $x$ for all $x\in \mathfrak{F}$.
\end{proof}

\begin{proof}[\scshape{Proof of Theorem \ref{thmP-4}} (Property P-4.).]
\emph{$h$-depth}. By Lemma \ref{lemmaHContinuous}, $D_{h}$ is continuous in $x$ so a fortiori, it is upper semicontinuous.

\emph{Random Tukey depth}. The case of $D_{RT}(y,P)\leq D_{RT}(x,P)$ is trivial. When $D_{RT}(y,P)> D_{RT}(x,P)$, the condition in \eqref{checkUSC} is
\begin{equation}\label{toVerifyRT}
\sup_{y:\; \|y-x\|<\delta} \min_{u\in\mathfrak{U}} D_{1}(\langle u,y\rangle ,P_{v})\leq \min_{u\in\mathfrak{U}} D_{1}(\langle u,x\rangle,P_{u}) + \epsilon.
\end{equation}
We verify the existence of a $\delta$ satisfying \eqref{toVerifyRT} for all 
\begin{equation}\label{eRT}
0<\epsilon\leq 1/2-D_{RT}(x,P).
\end{equation}
Note that if $D_{RT}(x,P)\geq 1/2,$ we are in the case of $D_{RT}(y,P)\leq D_{RT}(x,P)$. For the less interesting scenario in which $\epsilon>1/2-D_{RT}(x,P)$, the construction of $\delta$ satisfying \eqref{toVerifyRT} is more involved. Let $u\in\mathfrak{U}$ such that $D_{RT}(x,P)=D_{1}(\langle u,x\rangle,P_{u})$, and notice that $D_{RT}(y,P)\leq D_1(\langle u,y\rangle,P_u)$ for all $u\in\mathfrak{U}$. Additionally, $D_{RT}(y,P)>D_{RT}(x,P)=D_1(\langle u,x\rangle,P_u)$ implies $D_{1}(\langle u,y\rangle,P_{u})>D_{RT}(x,P)$. For $\epsilon$ satisfying (\ref{eRT}), $D_{RT}(x,P)=P_{u}(-\infty,\langle u,x\rangle]$ implies $D_{1}(\langle u,y\rangle,P_{u})=P_{u}(-\infty,\langle u,y\rangle]$ and, analogously,  $D_{RT}(x,P)=P_{u}[\langle u,x\rangle,-\infty)$ that $D_{1}(\langle u,y\rangle,P_{u})=[\langle u,y\rangle,-\infty).$
 With these observations, we see that \eqref{toVerifyRT} is achieved with $\delta<\sup\{\eta>0:P(B(\eta))\leq \epsilon\}$, where 
\[
B(\eta):=\{y\in\mathfrak{F}:D_{RT}(y,P)>D_{RT}(x,P)=D_1(\langle u,x\rangle,P_u), |\langle u,y-x\rangle|<\eta\}.
\]

\emph{Band depth and half-region depth}. \citet{Romo09} (Theorem 3) and \citet{LopezHalfRegion} (Proposition 6) prove that for all $x\in\mathfrak{F}$ and for all $\epsilon>0$, there exists a $\delta>0$ such that 
\[
\sup_{y:\;|\|y\|_{\infty}-\|x\|_{\infty}|<\delta}D(y,P)\leq D(x,P)+\epsilon
\]
for the respective depth constructions, $D=D_{J}$ and $D=D_{HR}$. Since $|\|y\|_{\infty}-\|x\|_{\infty}|\leq d(y,x)$ the proof is complete.

\emph{Modified band depth}. The case of $D_{MJ}(y,P)\leq D_{MJ}(x,P)$ is trivial. When $D_{MJ}(y,P)> D_{MJ}(x,P)$, the condition in \eqref{checkUSC} is
\begin{equation}\label{toVerifyMBD}
\sup_{y:\;\|x-y\|_{\infty}<\delta}\sum_{j=2}^{J}\mathbb{E}\bigl[\lambda\{v\in\mathcal{V}: y(v)\in [L_{j}(v),U_{j}(v)], x(v)\notin [L_{j}(v),U_{j}(v)]\}/\lambda(\mathcal{V})\bigr]\leq \epsilon,\quad J\geq 2.
\end{equation}
Taking $\delta<\sup\bigl\{\eta>0:\sum_{j=2}^{J}\mathbb{E}[\lambda\{v\in\mathcal{V}: x(v)\notin B_{j}(v), \min(|x(v)-L_{j}(v)|,|x(v)-U_{j}(v)|)<\eta\}]\leq\epsilon\lambda(\mathcal{V})\bigr\}$ ensures \eqref{toVerifyMBD} is satisfied.

\emph{Modified half region depth}. The case of $D_{MHR}(y,P)\leq D_{MHR}(x,P)$ is trivial. When $D_{MHR}(y,P)> D_{MHR}(x,P)$, the condition in \eqref{checkUSC} is
\begin{equation}\label{toVerifyMHR}
\sup_{y:\;\|x-y\|_{\infty}<\delta} \mathbb{E}[\lambda\{v\in\mathcal{V}: y(v)\leq X(v) \leq x(v)\}] \leq \epsilon.
\end{equation}
We verify the existence of a $\delta$ satisfying \ref{toVerifyMHR} for all $0<\epsilon\leq 1/2 - D(x,P)$. For the less interesting case of $\epsilon>1/2 - D(x,P)$, the construction of $\delta$ satisfying \ref{toVerifyMHR} is more involved. Let
\begin{eqnarray*}
\Gamma&:=& \{\eta>0:\mathbb{E}[\lambda\{v\in\mathcal{V}: (X(v)\leq x(v))\ind\{x\in\mathcal{A}\}, \\
& & \qquad \qquad (X(v)\geq x(v))\ind\{x\in\mathcal{B}\}, |x(v)-X(v)|]<\eta\}<\epsilon\lambda(\mathcal{V})\bigr\},
\end{eqnarray*}
where $\mathcal{A}:=\{x\in\mathfrak{F}:D(x,P)=\mathbb{E}[\lambda\{v\in\mathcal{V}:x(v)\leq X(v)\}]\}$ and $\mathcal{B}:=\{x\in\mathfrak{F}:D(x,P)=\mathbb{E}[\lambda\{v\in\mathcal{V}:x(v)\geq X(v)\}]\}$. Then taking $\delta<\sup\{\eta\in\Gamma\}$ ensures \eqref{toVerifyMHR} is satisfied.

\end{proof}

\begin{proof}[\scshape{Proof of Theorem \ref{thmP-5}} (Property P-5.).]

\emph{$h$-depth}. We obtain $D(f(x),P)>D(x,P)$ by simple calculation: $(\alpha(v))^2(x(v)-X(v))^2<(x(v)-X(v))^2$ for all $v\in L_\delta$ with $\lambda(L_\delta)>0$, hence
\[
D(f(x),P_{f(X)})=\frac{1}{h\sqrt{2\pi}}\EE\Bigl[\exp\Bigl\{-\frac{1}{2h^2}\bigl(\int_{L_\delta^c}(x(v)-X(v))^2dv+\int_{L_\delta}(\alpha(v))^2(x(v)-X(v))^2dv\bigr)\Bigr\}\Bigr]
\]
whilst $\displaystyle{
D(x,P)=\frac{1}{h\sqrt{2\pi}}\EE\Bigl[\exp\Bigl\{-\frac{1}{2h^2}\bigl(\int_{L_\delta^c}(x(v)-X(v))^2dv+\int_{L_\delta}(x(v)-X(v))^2dv\bigr)\Bigr\}\Bigr]}$.

\emph{Random Tukey depth}. The proof is by counterexample. Let $P$ be a discrete probability with $P[x_i]=1/3$ for $i=1,2,3$ and $x_1(v)>0,$  $x_2(v)=0$ and  $x_3(v)<0$ for all $v\in\mathcal{V},$ with $x_1$ and $x_3$ non-constant functions. Suppose for a contradiction that the following inequality is satisfied for $a=x_1$ and $a=x_3,$
\begin{eqnarray} \label{x1}D_{RT}(a,P_X)<D_{RT}(f(a),P_{f(X)}).\end{eqnarray}
If $a=x_1$ let's denote $b=x_3$ and else, if $a=x_3,$ $b=x_1.$
In general, as $\langle u, x_2\rangle=\langle u, f(x_2)\rangle=0,$
in order for the inequality (\ref{x1}) to be satisfied, any given $u\in\mathfrak{U}$ has to fulfil  either
$$\min\{0,\langle u, f(b)\rangle\}<\langle u ,f(a)\rangle<\max\{0,\langle u ,f(b)\rangle\} \mbox{ with }\langle u, f(b)\rangle\neq 0, \mbox{ or }$$
$$\langle u, f(a)\rangle=0\neq\langle u, f(b)\rangle, \mbox{ or }$$
\begin{eqnarray} \label{xp}\langle u, f(a)\rangle=\langle u, f(b)\rangle.\end{eqnarray}
However, in order for the inequality (\ref{x1}) to be simultaneously satisfied by $a=x_1$ and $a=x_3,$ only (\ref{xp}) can apply for each $u\in\mathfrak{U};$ but $\mu\{u: \langle u ,f(x_1)\rangle=\langle u , f(x_3)\rangle \}=0$ because, as $\alpha(v)>0$ for all $v\in\mathcal{V},$ $f(x_1)(v)>0$ and $f(x_3)(v)<0$ for all $v\in\mathcal{V}.$ Thus, (\ref{x1}) cannot be simultaneously satisfied by  $a=x_1$ and $a=x_3,$ which leads to contradiction.

\emph{Band depth, modified band depth, half region depth and modified half region depth}. The proof is by counterexample. We follow the counterexample of the random Tukey depth but state it here for the sake of completeness. Let $P$ be a discrete probability with $P[x_i]=1/3$ for $i=1,2,3$ and $x_1(v)>0,$  $x_2(v)=0$ and  $x_3(v)<0$ for all $v\in\mathcal{V},$ with $x_1$ and $x_3$ non-constant functions. As $\alpha(v)>0$ for all $v\in\mathcal{V},$ $f(x_1)(v)>0,$ $f(x_2)(v)=0$ and $f(x_3)(v)<0$ for all $v\in\mathcal{V}.$ In the case of the band depth and the modified band depth, for $j\in\{2,3\}$, the transformation simply shrinks the convex hull of any simplex over the $L_{\delta}$ region, whilst the probability of any simplex based on the transformation is the same as that of the original simplex to which it corresponds. It is thus immediate that $D(x_1,P_X)=D(f(x_1),P_{f(X)})$ for any $D\in\{D_{J},D_{MJ},D_{HR},D_{MHR}\}.$
\end{proof}

The proof of Theorem \ref{thmP-6} relies on the following definition.
\begin{definition}\label{defFrechet} \emph{[e.g.~\citet{vdVW1996}]}
For any map $\Psi:\mathbb{D}\mapsto \mathbb{K}$, with $\mathbb{D}$ and $\mathbb{K}$ normed spaces endowed with norms $\|\cdot\|_{\mathbb{D}}$ and $\|\cdot\|_{\KK}$ respectively, the Fr\'echet derivative of $\Psi$ \emph{(}if it exists\emph{)} is the linear continuous map $D\Psi_{a}:\mathbb{D}\mapsto \mathbb{K}$ such that 
\[
\left\|\Psi(a+b)-\Psi(a) - D\Psi_{a}(b)\right\|_{\KK} = o(\|b\|_{\mathbb{D}}).
\]
\end{definition}

\begin{proof}[\scshape{Proof of Theorem \ref{thmP-6} (Property P-6.)}.] \emph{$h$-depth}. Let $d_{\mathcal{P}}$ of Property P-6.~be the Prohorov metric or the bounded Lipschitz metric \citep[see e.g.][page 394]{DudleyRAP}. Then by the Portmanteau Theorem \citep[see e.g.][Theorem 11.3.3]{DudleyRAP}, 
$d_{\mathcal{P}}(Q,P)\rightarrow 0$ 
implies $|\int f(y)(P-Q)(dy)|\rightarrow 0$ for all $f\in BL (\mathfrak{F},d)$, where $BL(\mathfrak{F},d):=\left\{f:\mathfrak{F}\rightarrow \mathbb{R}:\; \|f\|_{BL}<\infty\right\}$, $\|f\|_{BL}=\|f\|_{L}+\|f\|_{\infty}$, and 
\[
\|f\|_{L}:=\sup_{z\neq y} \frac{|f(y)-f(z)|}{d(y,z)}.
\]
Given that, for any $x\in\mathfrak{F}$,
\[
|D_{h}(x,P)-D_{h}(x,Q)|=\Bigl|\int K_{h}(\|x-y\|_{\LL_{2}})P(dy) - \int K_{h}(\|x-y\|_{\LL_{2}})Q(dy)\Bigr|,
\]
it suffices by the previous observations to show that $K_{h}(\|x-\cdot\|_{\LL_{2}})\in BL(\mathfrak{F},d)$. First note
\[
\|K_{h}(\|x-\cdot\|_{\LL_{2}})\|_{\infty} = \sup_{y\in\mathfrak{F}}|K_{h}(\|x-y\|_{\LL_{2}})| = \sup_{y\in\mathfrak{F}}\Bigl|\frac{1}{h}\frac{1}{\sqrt{2\pi}}\exp\Bigl\{-\frac{\|x-y\|_{\LL_{2}}^{2}}{2h^{2}}\Bigr\}\Bigr|=(h\sqrt{2\pi})^{-1}<\infty.
\]
Thus it only remains to show $\|K_{h}(\|x-\cdot\|_{\LL_{2}})\|_{L}<\infty$. Taking $\Psi=K_{h}$, $a=x-z$, and $b=z-y$ in Definition \ref{defFrechet} yields
\[
\frac{\bigl| K_{h}(\|x-y\|_{\LL_{2}}) - K_{h}(\|x-z\|_{\LL_{2}}) - DK_{h,(x-z)}(z-y)\bigr|}{\|y-z\|_{\LL_{2}}}=o(1),
\]
Hence, to establish
\[
\sup_{z\neq y}\frac{\bigl| K_{h}(\|x-y\|_{\LL_{2}}) - K_{h}(\|x-z\|_{\LL_{2}})\bigr|}{\|y-z\|_{\LL_{2}}} < \infty,
\]
it is sufficient to show
\[
\sup_{z\neq y}\frac{|DK_{h,a}(z-y)|}{\|y-z\|_{_{\LL_{2}}}}<\infty.
\] 
Let $\psi(\cdot)=\|\cdot\|_{\LL_{2}}^{2}$ and $\varphi(\cdot)=\frac{1}{h}\frac{1}{\sqrt{2\pi}}\exp\{-\frac{(\cdot)}{2h^{2}}\}$. We can thus write $DK_{h,a}(z-y)=D_{a}(\varphi\circ\psi)(z-y)$, and by the chain rule of Fr\'echet derivatives, $D_{a}(\varphi\circ\psi)(b)=D_{a}\varphi((\psi)(b))\circ D_{a}\psi(b)$. We start by computing $D_{a}\psi(b)$. Setting $\Psi=\psi$ in Definition \ref{defFrechet} gives $\bigl|\langle a+b,a+b \rangle - \langle a,a \rangle - D\psi_{a}(b)\bigr| = o(\|b\|)$ and noticing that $\bigl|\langle a+b,a+b \rangle - \langle a,a \rangle - 2\langle a,b \rangle\bigr| = \langle b,b \rangle = \|b\|_{\LL_{2}}^{2} = o(\|b\|_{\LL_{2}})$, we conclude, $D\psi_{a}(b)=2\langle a,b \rangle = 2\langle x-z, z-y\rangle$. 

For an arbitrary $s\in\mathfrak{F}$, set $w=(\psi)(s)$, which belongs to $\mathbb{R}^{+}$, thus
\[
D_{a}\varphi(w)=-\frac{1}{2h^{3}}\frac{1}{\sqrt{2\pi}}\exp\Bigl\{-\frac{w}{2h^{2}}\Bigr\}.
\]
The chain rule delivers
\[
DK_{h,a}(z-y)=D_{a}(\varphi\circ\psi)(z-y)=-\frac{1}{h^{3}}\frac{1}{\sqrt{2\pi}}\exp\Bigl\{-\frac{\|y-z\|_{\LL_{2}}^{2}}{2h^{2}}\Bigr\}\langle x-z,z-y\rangle,
\]
hence
\begin{eqnarray*}
\sup_{z\neq y}\frac{|DK_{h,a}(z-y)|}{\|y-z\|_{\LL_{2}}} &=& \sup_{z\neq y}\frac{\Bigl|\frac{1}{h^{3}}\frac{1}{\sqrt{2\pi}}\exp\Bigl\{-\frac{\|z-y\|_{\LL_{2}}^{2}}{2h^{2}}\Bigr\}\langle x-z,z-y\rangle\Bigr|}{\|y-z\|_{\LL_{2}}}\\
&\leq & 
\sup_{z\neq y}\frac{1}{h^{3}}\frac{1}{\sqrt{2\pi}}\exp\Bigl\{-\frac{\|y-z\|_{\LL_{2}}^{2}}{2h^{2}}\Bigr\}\max{\{\|x\|_{\LL_{2}}^{2},\|y\|_{\LL_{2}}^{2},\|z\|_{\LL_{2}}^{2}\}} < \infty
\end{eqnarray*}
because $x,y,z\in\mathfrak{F}=\mathbb{L}_{2}$ implies they each have finite $\mathbb{L}_{2}$ norm.

\emph{Random Tukey depth}. $d_{\mathcal{P}}(Q,P)\rightarrow 0$ $P$-a.s.~for any metric $d_{\mathcal{P}}(\cdot,\cdot)$ metricising the topology of weak convergence, is equivalent to $Q\rightarrow P$ $P$-a.s., which in turn implies $Q_{u}\rightarrow P_{u}$ $P$-a.s.~for all $u\in\mathbb{H}$. 
 As $P$ is continuous and $u$ is drawn with a non-degenerate stationary Gaussian measure, $P_{u}$ is also continuous.
It follows that 
\[
\max\left\{|P_{u}(-\infty,\langle u,x \rangle]-Q_{u}(-\infty,\langle u,x \rangle]|,|P_{u}[\langle u,x \rangle,\infty)-Q_{u}[\langle u,x \rangle,\infty)|\right\} \rightarrow 0\;\; P-\text{a.s.},
\] and consequently, $|D_{1}(\langle u,x \rangle,P_{u})-D_{1}(\langle u,x \rangle,Q_{u})|\rightarrow 0$ $P$-a.s.~for any $u\in\mathbb{H}$. Then
\begin{eqnarray*}
|D_{RT}(x,P)-D_{RT}(x,Q)|&=&\bigl|\min_{u\in\mathfrak{U}}D_{1}(\langle u,x \rangle,P_{u})-\min_{u\in\mathfrak{U}}D_{1}(\langle u,x \rangle,Q_{u})\bigr|\\
												&\leq & \max_{u\in\mathfrak{U}}\bigl|D_{1}(\langle u,x \rangle,P_{u})-D_{1}(\langle u,x \rangle,Q_{u})\bigr| \rightarrow 0 \;\; P-\text{a.s.},
\end{eqnarray*}
where the inequality follows because, for any $w\in\mathfrak{U}$, $\min_{u\in\mathfrak{U}}D_{1}(\langle u,x \rangle, P_{u})\leq D_{1}(\langle w,x \rangle, P_{w})$, and likewise for $Q$.
The empirical case follows from the proof of Theorem 2.10 in \cite{randomTukey}.

\emph{Band depth}.  Since $d_{\mathcal{P}}(P,Q)$ metricises the weak topology, $d_{\mathcal{P}}(P,Q)<\delta \rightarrow 0$ is the same as writing $X_{\delta}\leadsto Y$ as $\delta\rightarrow 0$ where $\leadsto$ denotes weak convergence and $X_{\delta}$ and $Y$ are random variables $X_{\delta}:\Omega\rightarrow \mathfrak{F}$ and $Y:\Omega \rightarrow \mathfrak{F}$ such that, for any $A\in\mathcal{A}$, $P(A)=\mathbb{P}(X_{\delta}^{-1}(A))$ and $Q(A)=\mathbb{P}(Y^{-1}(A))$, where $\mathbb{P}$ is a probability on the underlying sample space $\Omega$. By the Portmanteau theorem \citep[e.g.][Theorem 11.3.3]{DudleyRAP}, $V_{N}\rightarrow_{d}V$ if and only if $\mathbb{E}f(V_{N})\rightarrow \mathbb{E}f(V)$ for all bounded Lipschitz functions $f$. Define $X_{\delta,1},\ldots,X_{\delta,J}$ to be i.i.d.~copies of $X_{\delta}$ and $Y_{1},\ldots,Y_{J}$ to be i.i.d.~copies of $Y$. Then, by the Portmanteau theorem, for any $\ell\in\{1,\ldots,j\}$ where $j\in\{2,\ldots,J\}$ and for any $(\alpha_{1},\ldots,\alpha_{j})\in\Delta^{j}$, since $f$ is bounded and continuous, there exists a $\delta<\delta_{\ell}$ such that
\[
\bigl|\mathbb{E}[f([\sum_{k\neq \ell}\alpha_{k}X_{\delta,k}] + \alpha_{\ell}X_{\delta,\ell})] - \mathbb{E}[f([\sum_{k\neq \ell}\alpha_{k}X_{\delta,k}] + \alpha_{\ell}Y_{\ell})]\bigr|<\delta/j.
\]
Hence
\begin{eqnarray*}
& & \bigl|\mathbb{E}[f(\sum_{k=1}^{j}\alpha_{k}X_{\delta,k})] - \mathbb{E}[f(\sum_{k=1}^{j}\alpha_{k}Y_{k})]\bigr| \\
& \leq & \sum_{\ell=1}^{j}\bigl|\mathbb{E}[f(\sum_{k\neq \ell}\alpha_{k}X_{\delta,k} + \alpha_{\ell}X_{\delta,\ell})] - \mathbb{E}[f(\sum_{k\neq \ell}\alpha_{k}X_{\delta,k} + \alpha_{\ell}Y_{\ell})]\bigr|<\delta.
\end{eqnarray*}
for all $\delta<\min\{\delta_{\ell}: \ell\in\{1,\ldots,j\}\}$. Letting $Z_{X(\delta),j}(\balpha):=\sum_{k=1}^{j}\alpha_{k}X_{\delta,k}$ and $Z_{Y,j}(\balpha):=\sum_{k=1}^{j}\alpha_{k}Y_{k}$, we conclude through a second application of the Portmanteau theorem that $Z_{X(\delta),j}(\balpha)\rightarrow_{d} Z_{Y,j}(\balpha)$ for any $j\in\{2,\ldots,J\}$ and any $\balpha\in\Delta^{j}$. Hence for every finite collection $\balpha_{1},\ldots,\balpha_{\ell}$ where $\balpha_{k}\in\Delta^{j}$ for each $k\in\{1,\ldots,\ell\}$, $\bigl(Z_{X(\delta),j}(\balpha_{1}),\ldots,Z_{X(\delta),j}(\balpha_{\ell})\bigr)\leadsto \bigl(Z_{Y,j}(\balpha_{1}),\ldots,Z_{Y,j}(\balpha_{\ell})\bigr)$. Here 
$\bigl(Z_{X(\delta),j}(\balpha_{1}),\ldots,Z_{X(\delta),j}(\balpha_{\ell})\bigr)$ is an arbitrary finite set of marginals (in the $\balpha$ index) $Z_{X(\delta),j}(\balpha):\Omega^{j}\rightarrow \mathfrak{F}$ of the stochastic process $Z_{X(\delta),j}:=\{Z_{X(\delta),j}(\balpha):\balpha\in\Delta^{j}\}$ which is the map $Z_{X(\delta),j}:\Omega^{j}\rightarrow \mathfrak{F}(\Delta^{j})=C(\mathcal{V},\Delta^{j})\subset \mathbb{L}^{\infty}(\mathcal{V}\times \Delta^{j})$, where $\mathbb{L}^{\infty}(\mathcal{V}\times \Delta^{j})$ is the space of bounded functions from $(\mathcal{V}\times \Delta^{j})$ to $\mathbb{R}$. Similarly, $\bigl(Z_{Y,j}(\balpha_{1}),\ldots,Z_{Y,j}(\balpha_{\ell})\bigr)$ is an arbitrary finite set of marginals of the stochastic process $Z_{Y,j}:=\{Z_{Y,j}(\balpha):\balpha\in\Delta^{j}\}$. Hence, in order to show that $Z_{X(\delta),j}\leadsto Z_{Y,j}$ for every $j\in\{2,\ldots,J\},$ it only remains by Theorem 1.5.4 of \citet{vdVW1996} to show that, for any $j\in\{2,\ldots,J\}$, $Z_{X(\delta),j}$ is asymptotically tight, i.e. for every $\xi>0$ there exists a compact set $K$ such that $\liminf_{\delta\rightarrow 0}P_{Z(\delta),j}\bigl(Z_{X(\delta),j}\in K^{\eta}\bigr)\leq 1-\xi$ for every $\eta>0$, where $P_{Z(\delta),j}$ is defined at every $A\in\mathcal{A}$ by $P_{Z(\delta),j}(A)=\mathbb{P}^{j}\bigl(Z_{X(\delta),j}^{-1}(A)\bigr)$. 

By Theorem 1.5.7 of \citet{vdVW1996}, $Z_{X(\delta),j}$ is asymptotically tight if and only if $Z_{X(\delta),j}(v,\balpha)$ is tight in $\RR$ for every $w=(v,\balpha)$, and there exists a semimetric $d_{w}$ on $\mathcal{W}=(\mathcal{V}\times \Delta^{j})$ such that ($\mathcal{W},d_{w}$) is totally bounded and $Z_{X(\delta),j}$ is uniformly $d_{w}$-equicontinuous in probability, i.e.~for every $\kappa,\varsigma>0$ there exists a $\gamma$ such that
\[
\limsup_{\delta\rightarrow 0} P_{Z(\delta),j}\left(\sup_{w,w':d_{w}(w,w')<\gamma}|Z_{X(\delta),j}(w)-Z_{X(\delta),j}(w')|>\kappa\right)<\varsigma.
\]
Tightness of $Z_{X(\delta),j}(v,\balpha)$ holds by completeness of $\mathfrak{F}$, which gives rise to tightness of $X_{\delta}$ and hence $Z_{X(\delta),j}$ because tightness is preserved under convex combinations. Since $\mathcal{V}$ is compact, so too is $\mathcal{W}$, hence $(\mathcal{W},d_{w})$ is totally bounded with respect to the $\ell_{1}$ norm. We have
\begin{eqnarray*}
&&Pr\left(\sup_{w,w':d_{w}(w,w')<\gamma}|Z_{X(\delta),j}(w)-Z_{X(\delta),j}(w')|>\kappa\right)\\ &\leq& Pr\left(\sup_{w,w':d_{w}(w,w')<\gamma}|Z_{X(\delta),j}(v,\balpha)-Z_{X(\delta),j}(v',\balpha)|>\kappa/2\right)\\ 
& & \qquad + Pr\left(\sup_{w,w':d_{w}(w,w')<\gamma}|Z_{X(\delta),j}(v',\balpha)-Z_{X(\delta),j}(v',\balpha')|>\kappa/2\right) = I+II.
\end{eqnarray*}
By the statement of Theorem \ref{thmP-6}, $\mathfrak{F}$ is the space of $d_{w}$-equicontinuous functions over $\mathcal{V}$. Since convex combinations of $d_{w}$-equicontinuous functions are $d_{w}$-equicontinuous, $Z_{X(\delta),j}(\cdot,\balpha)$ is $d_{w}$-equicontinuous with probability 1. It follows that 
for every $\kappa,\varsigma>0$, there exists a $\gamma>0$ such that $I<\varsigma/2$. Noting that $v'\in\mathcal{V}$ is fixed in II, taking $\gamma$ sufficiently small also gives rise to $II<\varsigma/2$, proving tightness. Asymptotic tightness is immediate because the bounds on $I$ and $II$ hold independently of $\delta$.

From here we know $Z_{X(\delta),j}\leadsto Z_{Y,j}$ for every $j\in\{2,\ldots,J\}$. It follows by Theorem 11.3.3 of \citet[][]{DudleyRAP} that there exists a $\eta(\delta)\searrow 0$ as $\delta \searrow 0$ such that $\rho(P_{Z(\delta),j},Q_{Z(Y),j})=M<\eta(\delta)$, where $Q_{Z(Y),j}(A)=\mathbb{P}^{j}\bigl(Z_{Y,j}^{-1}(A)\bigr)$, i.e.~for all $A\in\mathcal{A}$, $P_{Z(\delta),j}(A)\leq Q_{Z(Y),j}(A^{\xi})+\xi$ for all $\xi\in[M,\eta(\delta))$. Hence letting $B(x)=\cup\{A\in\mathcal{A}:x\in A\}$, we have $P_{Z(\delta),j}(B(x))\leq Q_{Z(Y),j}(B(x)^{\xi})+\xi$ for all $\xi\in[M,\eta(\delta))$ and by the symmetry of the Prohorov metric and the fact that $B(x)\subset B(x)^{\xi}$ for $\xi>0$ we conclude that $|P_{Z(\delta),j}(B(x))-Q_{Z(Y),j}(B(x))|\leq \xi<\eta(\delta)$. We have
\[
|D_{J}(x,P)-D_{J}(Q)|\leq \sum_{j=2}^{J}|P_{Z(\delta),j}(B(x))-Q_{Z(Y),j}(B(x))|<(J-1)\eta(\delta).
\]
Setting $\epsilon=(J-1)\eta(\delta)$, we see that the result follows by taking every $\delta_{\ell}$ in the above derivations equal to $\delta=\eta^{-1}(\epsilon/(J-1))$.

\emph{Modified band depth}. Let $[z(\cdot,\balpha):\balpha \in\Delta^{j}]$ be the set of all convex combinations of the elements of $\mathfrak{F}$ and $P_{Z(\delta),j}$ and $Q_{Y,j}$ be probability measures on that set, as defined in the proof of Theorem \ref{thmP-6} for the band depth. We have
\begin{eqnarray*}
\bigl| D_{MJ}(x,P)-D(x,Q)\bigr|&=&\bigl|\sum_{j=2}^{J}\frac{1}{\lambda(\mathcal{V})}\bigl(\EE\bigl[\lambda\{v\in\mathcal{V}:x\in \mathcal{S}_{j}(v,P)\}\bigr] -\EE\bigl[\lambda\{v\in\mathcal{V}: x\in \mathcal{S}_{j}(v,P)\}\bigr]\bigr)\Bigr| \\
													&\leq & \sum_{j=2}^{J}\Bigl|\frac{1}{\lambda(\mathcal{V})}\int \lambda \{v\in\mathcal{V}:x\in [z(v,\balpha):\balpha \in\Delta^{j}]\}(P_{Z(\delta),j}-Q_{Y,j})(dz)\Bigr|.
\end{eqnarray*}
But by compactness of $\mathcal{V}$, $\lambda \{v\in\mathcal{V}:x\in [z(v,\balpha):\balpha \in\Delta^{j}]\}$ is bounded and continuous in $z$ because $z\in \mathfrak{F}(\Delta^{j})=C(\mathcal{V}\times\Delta^{j})$. Hence $\bigl|D_{MJ}(x,P)-D(x,Q)\bigr|\rightarrow 0$ as $\delta\rightarrow 0$ by the Portmanteau Theorem \citep[][Theorem 11.3.3]{DudleyRAP} and the fact that $P_{Z(\delta),j}\rightarrow Q_{Y,j}$ as $\delta\rightarrow 0$, as demonstrated in the proof for the band depth.

\emph{Half region depth}. Take $d_{\mathcal{P}}(P,Q)=\rho(P,Q)$ where $\rho(P,Q)$ is defined as in the proof for the band depth. Suppose $\rho(P,Q)=M<\delta$ $P$-a.s., where $\delta>0$. Then for any $A\in\mathcal{A}$ and any $\eta\in[M,\delta)$, $P(A)-Q(A^{\eta})\leq \eta<\delta$. Let $E_{x}$ denote the epigraph of $x$ and let $H_{x}$ denote the hypograph of $x$. $|P(E_{x})-Q(E_{x})|\leq |P(E_{x})-Q(E_{x}^{\eta})|\leq \eta < \delta$ $P$-a.s.~and $|P(H_{x})-Q(H_{x})|\leq |P(H_{x})-Q(H_{x}^{\eta})|\leq \eta < \delta$ $P$-a.s., hence $\max\{|P(E_{x})-Q(E_{x})|,|P(H_{x})-Q(H_{x})|\}<\delta$ $P$-a.s. It follows that, for all $\epsilon>0$, $|D_{HR}(x,P)-D_{HR}(x,Q)|<\epsilon$ $P$-a.s.~as long as $d_{\mathcal{P}}(P,Q)<\delta$ $P$-a.s.~with $\delta=\epsilon$.

\emph{Modified half region depth}. Since $(\mathfrak{F},d)=(\mathcal{C}(\mathcal{V}),\|\cdot\|_{\infty})$ is separable and complete, $P$ and $Q$ are tight and by Theorem 11.3.5 and Corollary 11.6.4 of \citet{DudleyRAP}, $\rho(P,Q)=\alpha(X,Y)$ where $X$ and $Y$ are random variables with laws $P$ and $Q$ respectively, $\rho$ is the Prohorov metric defined and used throughout the proof of Theorem \ref{thmP-6} and $\alpha$ is the Ky-Fan metric, defined by $\alpha(X,Y):=\inf\{\eta>0:\Pr(d(X,Y)>\eta)\leq \eta\}$. Let $L$ be an arbitrary subset of $\mathcal{V}$ and let $X_{L}$ and $Y_{L}$ be the random variables $X$ and $Y$ defined over the restricted space with corresponding probability laws $P_{L}$ and $Q_{L}$ respectively. Since $P\rightarrow Q$, there exists a $\delta_{L}>0$ such that $\rho(P_{L},Q_{L})<\delta_{L}$, hence $\alpha(X_{L},Y_{L})<\delta_{L}$, hence $\Pr(d(X_{L},Y_{L})\geq \delta_{L})<\delta_{L}$ and for any Borel set $A_{L}$ of $C(L)$, if $X_{L}\in A_{L}$, then $Y_{L}\in A_{L}^{\delta_{L}}$, hence for any $L\subset \mathcal{V}$ and a sufficiently small $\delta_{L}$, 
$\{X_{L}(v)<x(v), Y_{L}(v)>x(v): v\in L\}$ and $\{X_{L}(v)>x(v), Y_{L}(v)<x(v): v\in L\}$ are events of probability zero under the joint law of $X_{L}$ and $Y_{L}$. By this argument,
\begin{eqnarray*}
|D_{MHR}(x,P)-D_{MHR}(x,Q)|&\leq& \max\Bigl\{ \bigl|\int \lambda\{v\in\mathcal{V}:y(v)\leq x(v)\}(P-Q)(dy)\bigr|, \\
													&&  \bigl| \int \lambda\{v\in\mathcal{V}:y(v)\geq x(v)\}(P-Q)(dy) \bigr| \Bigr\}
\end{eqnarray*}
with probability 1. Both terms in this expression converge to zero as $\delta\rightarrow 0$ by Theorem 11.3.3 of \citet{DudleyRAP} because $\lambda\{v\in\mathcal{V}:y(v)\leq x(v)\}$ and $\lambda\{v\in\mathcal{V}:y(v)\geq x(v)\}$ are continuous in $y$ and bounded by compactness of $\mathcal{V}$.
\end{proof}

% ##########################################################################################################################################################
% ##########################################################################################################################################################

\vspace{12pt}

\textbf{Acknowledgements:} Alicia Nieto-Reyes is grateful to the School of Mathematics at the 
University of Bristol for kind hospitality during the 2013/2014 academic 
year, during which time this work was carried out, and to the Spanish 
Ministerio de Ciencia e Innovaci\'on under grant MTM201128657-C0202 for 
financial support. Heather Battey is grateful for financial support from 
the EPSRC under grant EP/D063485/1. We thank Juan Cuesta-Albertos and 
Peter Green for their comments on a single-authored unpublished 
manuscript, from which the seeds of this paper originated, and Peter 
Hall for kindly reading the final draft. We also thank three anonymous referees for constructive suggestions.

\bibliographystyle{ims}
\bibliography{BiblioDepthProp}

\end{document}